\documentclass[a4paper,reqno]{amsart}

\usepackage{amsmath,amsfonts,amssymb}
\usepackage{graphics,graphicx,epsf}
\usepackage{xcolor,soul}
\usepackage{float}
\usepackage{hyperref}         
\usepackage{tikz}
\usepackage{enumitem}
\usepackage{tikz-cd}
\usetikzlibrary{arrows}

\newtheorem{theorem}{Theorem}[section]
\newtheorem{definition}[theorem]{Definition}
\newtheorem{proposition}[theorem]{Proposition}
\newtheorem{lemma}[theorem]{Lemma}
\newtheorem{corollary}[theorem]{Corollary}
\newtheorem{remark}[theorem]{Remark}

\def \R {\mathbb{R}}

\usepackage{geometry}
\numberwithin{equation}{section}

\addtocontents{toc}{\setcounter{tocdepth}{1}} 
\usepackage[symbol]{footmisc}

\title[]{Perturbation of parabolic equations with time-dependent linear operators: convergence of linear processes and solutions}

\author[M. Belluzi]{Maykel Belluzi$^\ast$}\thanks{$^\ast$Research 
	supported by FAPESP \# 2022/01439-5}

\address[M. Belluzi]{Universidade de São Paulo, Instituto de Ciências Matemáticas e de Computação, S\~{a}o
	Carlos SP, Brazil.}
\email{maykelbelluzi@icmc.usp.br}

\begin{document}
\maketitle

%%%%%%%%%%%%%%%%%%%%%%%%%%%%%%%%%%%%%%%%%%%%%%%%%%%%%
											%ABSTRACT
%%%%%%%%%%%%%%%%%%%%%%%%%%%%%%%%%%%%%%%%%%%%%%%%%%%%%
\vspace{-0.2cm}
\begin{abstract}

In this work we consider parabolic equations of the form
\[
(u_{\varepsilon})_t  +A_{\varepsilon}(t)u_{{\varepsilon}} = F_{\varepsilon} (t,u_{{\varepsilon} }), 
\]
where $\varepsilon$ is a parameter in $[0,\varepsilon_0)$ and $\{A_{\varepsilon}(t), \ t\in \mathbb{R}\}$ is a family of uniformly sectorial operators. As $\varepsilon \rightarrow 0^{+}$, we assume that the equation converges to
\[
u_t  +A_{0}(t)u_{} = F_{0} (t,u_{}).
\]
{The time-dependence found on the linear operators $A_{\varepsilon}(t)$ implies that linear process is the central object to obtain solutions via variation of constants formula}. Under suitable conditions on the family $A_{\varepsilon}(t)$ and on its convergence to $A_0(t)$ when $\varepsilon \rightarrow 0^{+}$, we obtain a Trotter-Kato type Approximation Theorem for the linear process $U_{\varepsilon}(t,\tau)$ associated to $A_{\varepsilon}(t)$, estimating its convergence  to the linear process $U_0(t,\tau)$ associated to $A_0(t)$. Through the variation of constants formula and assuming that $F_{\varepsilon}$ converges to $F_0$, we analyze how this linear process convergence is transferred to the solution of the semilinear equation. We illustrate the ideas in two examples. First a reaction-diffusion equation in a bounded smooth domain $\Omega \subset \mathbb{R}^{3}$
\[\begin{split}
& (u_{\varepsilon})_t - div (a_{\varepsilon} (t,x) \nabla u_{\varepsilon}) +u_{\varepsilon} = f_{\varepsilon} (t,u_{\varepsilon}), \quad x\in \Omega,  t> \tau,
\\
%& \partial_n u_{\varepsilon} =0,  \hspace{4.45cm} x\in \partial \Omega,
\end{split}
\] 
where $a_\varepsilon$ converges to a function $a_0$, $f_{\varepsilon}$ converges to $f_0$. 
%and Neumann boundary condition is assumed. 
We apply the abstract theory in this example, obtaining convergence of the linear process and solution. As a consequence, we also  obtain upper-semicontinuity of the family of pullback attractors associated to each problem. {The  second example is a nonautonomous strongly damped  wave equation %in a bounded smooth domain $\Omega \subset \R^{n}$,
\[
u_{tt}+(-a(t) \Delta_D) u + 2 (-a(t)\Delta_D)^{\frac12} u_t = f(t,u),  \quad 
x\in \Omega,
 t>\tau
,\]
where $\Delta_D$ is the Laplacian operator with Dirichlet boundary conditions in a domain $\Omega$
and we analyze convergence of solution as we perturb the fractional powers of the associated linear operator.
}

%MSC

\end{abstract}

\tableofcontents

\allowdisplaybreaks

%%%%%%%%%%%%%%%%%%%%%%%%%%%%%%%%%%%
				%SECTION 1 - INTRODUCTION
%%%%%%%%%%%%%%%%%%%%%%%%%%%%%%%%%%%

\section{Introduction}\label{Introduction}
	
	In the present paper we study \textit{singularly nonautonomous} semilinear parabolic problems of the form
	\begin{equation}\label{Pe}
	\begin{split}
		&(u_{\varepsilon})_t  +A_{\varepsilon}(t)u_{{\varepsilon}} = F_{\varepsilon} (t,u_{{\varepsilon} }), \quad { } t>\tau, \\
		&u_{\varepsilon}(\tau) = u^{\tau} \in Y \stackrel{}{\hookrightarrow} X,
	\end{split}
	\end{equation}
	where $\varepsilon \in [0,\varepsilon_0)$ is a parameter, $X$ is a Banach space, $A_{\varepsilon}(t): D(A_{\varepsilon}(t)) \subset X \to X$ is a family of sectorial operators (with a certain uniformity in $t$ that we shall specify latter), $Y$ is a Banach Space continuously embedded in  $X$, which we denote by $Y \stackrel{}{\hookrightarrow} X$, and $F_{\varepsilon}: \mathbb{R} \times Y \to X$ is a nonlinearity.
	
	The term \textit{singularly nonautonomous} expresses the fact that this family $A_{\varepsilon} (t)$ is time-dependent, as a counterpart to the  semilinear problems where $A_{\varepsilon}(t) = A_{\varepsilon}$, which we refer as \emph{nonsingular}. {This terminology, adopted for instance in \cite{BCNS_2,CN}, is not unanimous and, in the case we are considering here, does not refer to any discontinuity or blow-up in time, which  can mean in other contexts.} We shall adopt it in order to easily distinguish between the case studied to the well established case where there is no time-dependence on the linear operators. 
	
	As $\varepsilon \rightarrow 0^{+}$, Problem \eqref{Pe} approaches to what we refer as \emph{limiting problem}
	\begin{equation*}\label{P0}
	\begin{split}
	&(u_0)_t  +A_{0}(t)u_{0} = F_{0} (t,u_{0}), \quad t>\tau,\\
	&u_0(\tau) = u^{\tau} \in Y,
	\end{split}
	\end{equation*}
	whose solution is denoted by $u_0(t) = u_0(t,\tau,u^{\tau})$ and referred as \emph{limiting solution}.

	For each $\varepsilon \in [0,\varepsilon_0)$ and under suitable conditions on the family $\{ A_{\varepsilon}(t), t\in \mathbb{R} \}$ and on the nonlinearity $F_{\varepsilon}$, Problems \eqref{Pe} are well-posed. 
	 We are interested in investigating the behavior of the solution $u_{\varepsilon}(t)$ of \eqref{Pe} as $\varepsilon \rightarrow 0^{+}$, comparing it to the limiting solution $u_0(t)$ and providing a rate of convergence for those solutions in terms of $\varepsilon$. In order to obtain this convergence, we must first study the associated linear problem.

	This type of analysis has already been done when the linear operators in \eqref{Pe} do not depend on time, that is, $A_{\varepsilon} (t) \equiv A_{\varepsilon}$, $\varepsilon \in [0,\varepsilon_0)$. In this case, each operator $-A_{\varepsilon}$ generates a linear {semigroup}, $\{T_{-A_{\varepsilon}}(t) \in \mathcal{L}(X), t\geq 0\}$, that plays an essential role in solving the semilinear problem. Under suitable assumptions on $F_{\varepsilon}$, the \textit{nonsingular} problem 
	\begin{equation}\label{not_singular}
	(u_{\varepsilon})_t+A_{\varepsilon}u_{\varepsilon}=F_{\varepsilon}(t,u_{\varepsilon}), \ t>\tau, \quad u_{\varepsilon}(\tau)=u^{\tau} \in Y, \quad \varepsilon \in [0,\varepsilon_0),
	\end{equation}  is locally solved by
	\begin{equation}\label{FVC_semigrupo}
	u_{\varepsilon}(t) = T_{-A_{\varepsilon}}(t-\tau)u^{\tau} + \int_{\tau}^{t} T_{-A_{\varepsilon}}(t-s) F_{\varepsilon}(s,u_{\varepsilon}(s))ds
	\end{equation}
	and we refer to the above expression as \textit{variation of constants formula}.

	In papers such as \cite{Arrieta_1999, ArrietaBezerraCarvalho_2013, ArrietaCarvalho_2004, Carbone_et_al_2008, CarvalhoPires_2017, CarvalhoPiskarev_2004, AbreuCarvalho_2004} a general routine was conceived and applied in order to guarantee convergence of solutions of Problems \eqref{not_singular} as $\varepsilon \to 0^{+}$. This routine is based in a detailed study of the behavior of the linear part. Precisely, the routine  consists in first studying the convergence of the linear operator $A_{\varepsilon}^{-1}$ to $A^{-1}_0$. This information is then used to obtain the convergence of the resolvent operator $(\lambda+A_{\varepsilon})^{-1}$ to $(\lambda+A_0)^{-1}$ in some sector and from the resolvent convergence one obtains convergence of the linear semigroup $T_{-A_{\varepsilon} }(\cdot)$ to $T_{-A_0}(\cdot)$. By using the variation of constants formula \eqref{FVC_semigrupo}, one can also prove the convergence of the solutions to the limiting solution as a consequence of the linear semigroup convergence.

	If the equation $u^{\varepsilon}_t+A_{\varepsilon}u^{\varepsilon}=F_{\varepsilon}(u^{\varepsilon})$ is autonomous ($F_{\varepsilon}$ does not depend on time), one can continue the analysis and derive the upper semi-continuity of the family of global attractors $\{ \mathcal{A_{\varepsilon}} \}_{\varepsilon \in [0,1]}$ and even lower semi-continuity under suitable structural hypothesis on the limiting attractor $\mathcal{A}_0$. This is done for instance in  \cite{ArrietaBezerraCarvalho_2013}.

	A careful analysis of those papers allows us to conclude that a huge effort goes in the direction of ensuring that the linear semigroup $T_{-A_{\varepsilon}}(\cdot)$ converges to $T_{-A_0}(\cdot)$ in an appropriate sense. From this type of Trotter-Kato Approximation Theorem, the convergence of the other elements being studied follows.

	The situation changes  when we consider singularly nonautonomous problem, since the linear semigroup is not the central element in obtaining the solution of the semilinear problem, as we discuss next. However, we are still able to elaborate for the singularly nonautonomous case (Problems \eqref{Pe}) a routine similar to the one mentioned in the articles above to treat the matter of convergence for the problems. An approach like this one for the singular nonautonomous case does not exist in the literature so far, and with the results we present in this paper, we shall be able to study perturbation of singularly nonautonomous problem, incorporating several different examples that appears in the literature.

	The only matter that we shall not address in this paper is the lower semi-continuity of the pullback attractor associated to Problems \eqref{Pe}, whenever we are able to prove that they exist. We do not pursuit this result due to the fact that there might not exist an elliptic associated problem for the limiting equation 
	\[
	(u_0)_t+A_0(t)(u_0) = F_0(t,u_0)
	\]
	and we usually can not derive information on the structure of the pullback attractor unless we require some simplifying assumptions on $A_0(t)$ and $F_0(t,\cdot)$ with respect to the time-dependence. Since this is not the purpose of this article, we do not look for a result on  lower semi-continuity of the family of pullback attractors. Nevertheless, upper-semicontinuity of attractors will be obtained as a consequence of the convergences established for the solutions.
	
	The main difference between the case where $A_{\varepsilon} (t) \equiv A_{\varepsilon}$   to the singularly nonautonomous comes from the fact that instead of a linear semigroup associated to $-A_{\varepsilon}$ that provides solutions for the semilinear problem through the variation of constant formula \eqref{FVC_semigrupo}, we will have a two parameter family of linear operators 
	$$\{U_{\varepsilon} (t,\tau) \in \mathcal{L}(X),\  t\geq\tau, \ \tau\in \mathbb{R} \}$$ that will be essential in describing the solution for the semilinear problem. The existence of such family associated to $\{A_{\varepsilon} (t), t\in \mathbb{R}\}$ was established almost simultaneously by Sobolevski\u{\i} \cite{sobol} and Tanabe \cite{Tanabe_1959, Tanabe_1960, Tanabe_1960_2}. 
	This family $U_{\varepsilon}(t,\tau)$ has properties similar to the ones presented by the linear semigroup in the autonomous case. In particular, there is an equivalent variation of constant formula that provides solutions for \eqref{Pe}  given by
	\begin{equation}\label{FVC_processo}
	u_{\varepsilon}(t) = U_{\varepsilon}(t,\tau)u_{\varepsilon}^{\tau} + \int_{\tau}^{t} U_{\varepsilon}(t,s) F_{\varepsilon}(s,u_{\varepsilon}(s))ds.
	\end{equation}
	
	Taking this into account, the outline we adopt to treat perturbation of singularly nonautonomous parabolic problems consists in the following steps:
	\begin{enumerate}[label={{(\roman*)}},ref=(\roman*)]
		\item \label{(I)} First we prove that, for each fixed time $t \in \mathbb{R}$, the linear operator $A_{\varepsilon}(t)^{-1}$ converges in an appropriate sense to the linear operator $A_0(t)^{-1}$.
		We also establish the rate of such  convergence in terms of specific characteristics of the problem.
		
		\item We use the previous information to obtain the resolvent convergence of  $(\lambda+A_{\varepsilon}(t))^{-1}$ to $(\lambda + A_{0}(t))^{-1}$ in a sector common to all the resolvent sets of all linear operators.
	
		\item Through a well-known formulation for analytic semigroups in terms of its resolvent, we transfer the resolvent convergence to the linear semigroup generated by $-A_{\varepsilon}(t)$, for a fixed $t\in \mathbb{R}$, that is, we obtain the convergence (with rate) of $T_{-A_{\varepsilon}(t)} (\cdot)$ to $T_{-A_0(t)} (\cdot)$.
		
		\item Using the formulations of the linear process $U_{\varepsilon}(t,\tau)$ in terms of $A_{\varepsilon}(t)$ and $T_{-A_{\varepsilon}(t)} (\cdot)$ (developed in \cite{sobol} which we discuss in the sequel), we obtain the convergence (with rate) of $U_{\varepsilon} (t,\tau)$ to $U_0(t,\tau)$. This result is presented  Theorem \ref{T_Conv_Proc}.

		\item \label{(V)} Using the variation of constants formula \eqref{FVC_processo}, we obtain in Theorem \ref{T_Conv_solution} the convergence (with rate) of the solution  $u_{\varepsilon} (\cdot)$ to the solution $u_0(\cdot)$.

		\end{enumerate}

	To attend the program proposed, 
	this paper is structured in the following manner: In Section \ref{S:Setting_and_main_results} we present the assumptions required for the family of linear operators $\{A_{\varepsilon}(t), t\in \mathbb{R}\}$ and for the nonlinearities $F_{\varepsilon}$ that allow us to prove the results on convergence. We also enunciate in this section the main abstract results on convergence: Theorem \ref{T_Conv_Proc} on the convergence of the linear process and Theorem \ref{T_Conv_solution} on the convergence of the solutions of \eqref{Pe} as $\varepsilon \to 0^{+}$.  Their proofs are postponed to Section \ref{S_Estimates_rates_convergence} and they depend on following steps \ref{(I)} to \ref{(V)} mentioned above.  We then apply those results in two different examples. First, in Section \ref{S:Case_I_application_I}, we consider a family of reaction-diffusion equations in a fixed bounded smooth domain $\Omega \subset \mathbb{R}^{3}$
	 \[\begin{split}
	 & (u_{\varepsilon})_t - div (a_{\varepsilon} (t,x) \nabla u_{\varepsilon}) +u_{\varepsilon} = f_{\varepsilon} (t,u_{\varepsilon}), \quad x\in \Omega,  \ t> \tau,\\
	 & \partial_n u_{\varepsilon} =0,  \hspace{5.11cm} x\in \partial \Omega.
	 \end{split}
	 \] 
	Assuming that $a_\varepsilon$ converges to $a_0$ and $f_{\varepsilon}$ converges to $f_0$, we derive in this section all the abstract conditions required in Theorems \ref{T_Conv_Proc} and \ref{T_Conv_solution} that ensures convergence of the solution $u_{\varepsilon}$ to $u_0$, as $\varepsilon \to 0^{+}$. Moreover, under an additional dissipation assumption on the nonlinearities $f_{\varepsilon}$, we prove that each problem is globally well-posed, defines a nonlinear dissipative process with pullback attractor $\{\mathcal{A}_{\varepsilon} (t) \subset Y, \ t\in \mathbb{R}\}$ and we prove this family of pullback attractors is upper-continuous in $\varepsilon=0$.
	Finally, in Section \ref{S_application_II}, we apply the abstract theory to a nonautonomous strongly damped wave equation and its fractional approximations in the sense of \cite{BezerraCarvalhoNascimento2020}.

	\

	Before we proceed to the goals proposed, we mention two points that are important to take into account. The first one concerns the linear operators $A_{\varepsilon}(t)$. For this paper, we shall consider a situation where  the domain $D(A_{\varepsilon}(t)) = D_{\varepsilon}$ remains fixed in $t$ and  the phase space $X$ where the linear operator is defined  remains fixed in $t$ and $\varepsilon$. This assumption holds for several problems, as we shall see in applications. The situation where the domain of the linear operator is time-dependent and the phase space changes with  $\varepsilon$ or time shall be addressed in future works.

	The second point that we want to highlight is the motivation behind considering \textit{singularly nonautonomous} problems. In general, an evolution system in a Banach space $X$ can be represented by an equation
	\begin{equation}\label{Prob_nao_linear}
	u_t = f(t,u), \quad t> \tau, \qquad u(\tau)=u^{\tau} \in Y,
	\end{equation}
	
	\vspace{0.1cm}\noindent where  $Y \hookrightarrow X$ and  $f: \mathcal{U} \subset \mathbb{R} \times Y \to X$. However, the function $f$ can be highly nonlinear, which makes it difficult to study the problem. To simplify it, we can approximate the above equation around a state $u_0$ by a linear (or semilinear) evolution equation, and then use the several tools mentioned above (and others in the existent literature) to treat semilinear problems. 
	
	This approximation is obtained by considering the Taylor polynomial of $f$ around the state $u_0$ (assuming that $f$ has the necessary regularity), that is,
	\begin{equation*}
	f(t,u_0+z) = f(t,u_0) + \frac{ \partial f}{\partial u} (t,u_0) z + g(t,z),
	\end{equation*}
	where $g(t,z) = o(\|z\|_{Y})$ when $\|z\|_Y \to 0$ and $\frac{ \partial f}{\partial u} (t,u_0) \in \mathcal{L}(Y,X)$ is the Frechet Derivative of $f$ with respect to the second variable. If we denote $z(t) = u(t)-u_0$ and $-A(t) =  \frac{ \partial f}{\partial u} (t,u_0)$, Problem \eqref{Prob_nao_linear} becomes 
	\begin{equation*}\label{Prob_apos_Taylor}
	z_t + A(t)z = f(t,u_0) + g(t,z), \quad t> \tau, \qquad z(\tau)=0,
	\end{equation*}
	which is \textit{singularly nonautonomous} and is in the same format as Problem \eqref{Pe}. Therefore, singularly nonautonomous evolution equations seems to  be a good tool to model several real life phenomena and compels the efforts in the direction of describing its dynamics.

\vspace{0.1cm}

%%%%%%%%%%%%%%%%%%%%%%%%%%%%%%%%%%%%%%%%%%%%%%%%%%%%
%SECTION 2
%%%%%%%%%%%%%%%%%%%%%%%%%%%%%%%%%%%%%%%%%%%%%%%%%%%%

\section{Functional setting and main results}\label{S:Setting_and_main_results}

In the sequel we provide conditions  on the family of linear operators $\{A_{\varepsilon}(t), \ t\in \mathbb{R} \}$ that ensure existence of the linear process $\{U_{\varepsilon}(t,\tau) \in  \mathcal{L}(X), \ t\geq \tau, \ \tau\in \mathbb{R}\}$ as well as convergence of  $U_{\varepsilon}(t,\tau)$ to $U_0(t,\tau)$ as $\varepsilon \to 0^{+}$. 
Once we have convergence of the linear parts established, we provide conditions for the nonlinearities $F_{\varepsilon}$ that guarantee  convergence of the solutions of \eqref{Pe} as $\varepsilon \to 0$.

\subsection{A type of Trotter-Kato Approximation Theorem  for the Linear processes $U_{\varepsilon}(t,\tau)$}
Consider the abstract \textit{singlularly nonautonomous} semilinear  problem \eqref{Pe} and assume that, for each  $\varepsilon \in [0,\varepsilon_0)$, $\{A_{\varepsilon} (t), \ t\in \R\}$ is a family of linear operators in  $X$ satisfying: 

\vspace{0.3cm}

\begin{enumerate}[label={\textbf{(P.\arabic*)}},ref=(P.\arabic*)]
	
	\item \label{P_1} The operator $A_{\varepsilon}(t):D (A_{\varepsilon}(t)) \subset X_{}\to X_{}$ is a closed densely defined linear operator, the domain $D_{\varepsilon}=D(A_{\varepsilon}(t))$ is fixed  in time (but it can change with $\varepsilon$) and there are 
	constants ${C}>0$ and $\mathbf{\varphi} \in (\frac{\pi}{2}, \pi)$ (independent of $\varepsilon \in [0,\varepsilon_0)$ and $t\in \R$) such that 
			\begin{equation*}
			\Sigma_{\varphi} \cup \{0\} \subset \rho (-A_{\varepsilon}(t)), \quad \mbox{ for all } \varepsilon \in [0,\varepsilon_0) \mbox{ and } t\in \R,
			\end{equation*}
			
			\noindent where $\Sigma_{\varphi} = \left\{  \lambda \in \mathbb{C}; |\arg \lambda| \leq \varphi \right\}$ and
			\begin{equation}\label{res_est_X}
			\| (\lambda I + A_{\varepsilon}(t) )^{-1}\|_{\mathcal{L}(X)} \leq \frac{C}{|\lambda|^{} +1}, \quad \mbox{ for all } \lambda\in \Sigma_{\varphi} \cup\{0\}.
			\end{equation}
	
			We say in this case that the family $A_{\varepsilon}(t)$ is \emph{uniformly sectorial}.

	\item \label{P_2} The operators $A_{\varepsilon}(t)$ have the following  regularizing property: its resolvent has its image on the Banach space $Y \hookrightarrow X$ and there exists $\beta\in (0,1]$ such that 
	\begin{equation}\label{res_est_Y}
	\| (\lambda I + A_{\varepsilon}(t) )^{-1}\|_{\mathcal{L}(Y)} \leq \frac{C}{|\lambda|^{} +1}, \quad \mbox{ for all } \lambda\in \Sigma_{\varphi} \cup\{0\},
	\end{equation}
	and
	\begin{equation}\label{res_est_X,Y}
	\| (\lambda I + A_{\varepsilon}(t) )^{-1}\|_{\mathcal{L}(X,Y)} \leq \frac{C}{|\lambda|^{\beta} +1}, \quad \mbox{ for all } \lambda\in \Sigma_{\varphi} \cup\{0\}.
	\end{equation}

	\item \label{P_3} There are constants $C>0$ and $0<\delta \leq1$ (independent of $\varepsilon \in [0,\varepsilon_0)$) such that, for any
	$t,\tau,s \in {\mathbb R}$,
	\begin{equation*}\label{1.16}
	\|[A_{\varepsilon}(t)-A_{\varepsilon}(\tau)]A_{\varepsilon}(s)^{-1}\|_{\mathcal{L}(X)} \leq C|t-\tau|^{\delta}.
	\end{equation*}
	We say that the function $\R \ni t\mapsto A_{\varepsilon}(t)A_{\varepsilon}(s)^{-1}  \in \mathcal{L}(X)$ is $\delta-$\emph{uniformly H\"{o}lder continuous}.
	
\end{enumerate}

\

	Conditions \ref{P_1} and \ref{P_2} state that each operator $A_{\varepsilon} (t)$, $\varepsilon \in [0,\varepsilon_0)$ and $t\in \mathbb{R}$, is sectorial and we can guarantee the existence of a common sector in the spectrum of them all as well as uniform estimate in this sector. Those properties can be seen as a uniform parabolicity for the family of  linear operators. Moreover, to say that the resolvent of $A(t)$ has its image in $Y$ means that $D(A(t)) \subset Y$, since $(\lambda+A(t))^{-1}:X\to D(A(t))$.

	Condition \ref{P_3} states that the H\"{o}lder exponent for the maps $t\mapsto A_{\varepsilon}(t)A_{\varepsilon} (s)^{-1} \in \mathcal{L}(X)$ can be chosen uniformly among all families and, as a consequence of this property, 
	$$\| A_{\varepsilon}(t)A_{\varepsilon}(\tau)^{-1} \|_{\mathcal{L} (X)  } \leq C, \quad  \mbox{ for all } (t,\tau) \mbox{ in a compact set and } \varepsilon\in [0,\varepsilon_0).
	$$ 
	
	In this case, the graph norms defined by the operators $A_{\varepsilon}(t)$ and $A_{\varepsilon}(\tau)$ in $D_{\varepsilon}$,
	 $$\| \cdot \|_{D(A_{\varepsilon}(t))} = \| A_{\varepsilon}(t) \cdot \|_{X_{}} \quad \mbox{ and } \quad \| \cdot \|_{D(A_{\varepsilon}(\tau))} = \| A_{\varepsilon}(\tau) \cdot \|_{X_{}},$$
	 respectively, are equivalent. We shall refer to both norms as $\| \cdot \|_{X_{}^{1}}$.
	
	From conditions \ref{P_1} to \ref{P_3} will be able to derive uniform estimates in $\varepsilon$ for the semigroups and linear process associated to the family 
	$\{A_{\varepsilon}(t), \ t\in \mathbb{R}\}$. Nevertheless, in order to obtain properties of convergence as we make $\varepsilon \to 0^{+}$, we shall require further conditions on the linear operators that connect the different problems being studied. Those conditions are stated next:

	\begin{enumerate}[label=\textbf{(P.\arabic*)},ref=(P.\arabic*)]
		\setcounter{enumi}{3} 
		
		\item \label{P_4} There exists a continuous function $\xi:[0,\varepsilon_0) \to \mathbb{R}^{+}$ with $\xi (0)=0$ such that  
		\begin{equation*}
		\sup_{t, \tau \in \mathbb{R}}\| A_{\varepsilon}(t)A_{\varepsilon}(\tau)^{-1} - A_{0}(t)A_0(\tau)^{-1}\|_{\mathcal{L}(X)} \leq \xi (\varepsilon).
		\end{equation*}

		\item \label{P_5} There exists a continuous function $\eta: [0,\varepsilon_0) \to \mathbb{R}^{+}$ with $\eta(0)=0$ such that 
		\begin{equation*}
		\sup_{t \in \mathbb{R}} \| A_{\varepsilon} (t)^{-1} -  A_0(t)^{-1} \|_{\mathcal{L} (X,Y)} \leq  \eta (\varepsilon).
		\end{equation*}
		%We refer to this $\eta$ as the \textit{rate of convergence of the linear operators}.
		
	\end{enumerate}

\

For a fixed $\varepsilon \in [0,\varepsilon_0)$ and $\tau \in \R$, each operator $A_{\varepsilon} (\tau)$ is sectorial  with $\Sigma_{\varphi}\cup\{0\}$ in the resolvent of $-A_{\varepsilon}(\tau)$. Henceforth, $-A_{\varepsilon} (\tau)$ generates an analytic semigroup which we denote by $T_{-A_{\varepsilon} (\tau)} (\cdot)$ (see \cite[Theorem 1.5.2]{pazy}) given by

\begin{equation}\label{semigrupo}
T_{-A_{\varepsilon}(\tau)}(t) = \dfrac{1}{2\pi i} \int_{\Gamma} e^{\lambda t} (\lambda + A_{\varepsilon}(\tau))^{-1}  d\lambda,
\end{equation}
\vspace{0.1cm}

\noindent where $\Gamma$ is the contour of $\Sigma_{\varphi}$  and it is oriented with increasing imaginary part. This linear semigroup solves the linear and homogeneous differential equation
\[
(u_{\varepsilon})_t +A_{\varepsilon}(\tau)  u_{\varepsilon} = 0, \ t>0, \quad u_{\varepsilon}(0)=u^{0} %\in Y
,
\]
by considering $u_{\varepsilon}(t,u^{0})= T_{-A_{\varepsilon}(\tau)} (t) u^{0}.$

However, the family $\{T_{-A_{\varepsilon}(\tau)} (t) \in \mathcal{L}(X), \ t \geq 0\}$  is not enough to describe the evolution of the system associated to \eqref{Pe}. We must obtain a two parameter family $\{U_{\varepsilon}(t,\tau ) \in \mathcal{L}(X), t\geq \tau\}$ of linear operators associated to $A_{\varepsilon}(t)$ that plays in the singularly nonautonomous case 
a similar role as the semigroup in the nonsingular case, that is, we should expect $U_{\varepsilon} (t,\tau)$ to recover the solution of  the homogeneous equation
\begin{equation}\label{homogeneous_diff_eq}
(u_{\varepsilon})_t+A_{\varepsilon}(t)u_{\varepsilon} =0, \mbox{ } t>\tau;\quad u_{\varepsilon}(\tau) = u^{\tau},%\in Y,
\end{equation}
by considering $u_{\varepsilon }(t) = U_{\varepsilon} (t,\tau)u^{\tau}$. In other words, we expect that $\partial_t U_{\varepsilon} (t,\tau) = -A_{\varepsilon } (t) U_{\varepsilon}(t,\tau)$. As a matter of fact, we search for the existence of a family of linear operators $\{ U_{\varepsilon} (t,\tau)\in \mathcal{L}(X), \ t\geq \tau \}$ with the following properties:

\begin{definition}\label{evol_process}
	Let $X$ be a Banach space. A family $\{U_{\varepsilon}(t,\tau) \in \mathcal{L}(X), t\geq \tau\}$ of bounded linear operators is a \emph{linear process associated to $A_{\varepsilon}(t):D_{\varepsilon} \subset X \to X$} if 
	
	\begin{enumerate}
		\item $U_\varepsilon(t,t) = I$ and $U_{\varepsilon}(t, s)U_{\varepsilon}(s,\tau) = U_{\varepsilon}(t,\tau)$, for all $\tau\leq s \leq t$.
	
		\item $(t,\tau,x) \mapsto U_{\varepsilon}(t,\tau)x$ is continuous for $t\geq \tau$ and for all $ x \in X$. 
		
		\item There exist $C, T >0$ such that $\| U_{\varepsilon}(t,\tau)\|_{\mathcal{L}(X)} \leq C$, for all $0\leq t-\tau \leq T$.
		
		\item $U_{\varepsilon}(t,\tau):X \rightarrow D_{\varepsilon}$ and $(\tau, \infty) \ni t \mapsto U_{\varepsilon}(t,\tau)x \in X$ is differentiable for each $x\in X$.
		
		\item The derivative $\partial_t U_{\varepsilon}(t,\tau) $ is a bounded linear operator in $X$, 
		$$\partial_t U_{\varepsilon}(t,\tau) =-A_{\varepsilon}(t) U_{\varepsilon}(t,\tau)$$ 
		and, for  $T>0$, there exists $C=C(T)>0$ such that 
		$$\| \partial_t U_{\varepsilon}(t,\tau) \|_{\mathcal{L}(X)} \leq C(t-\tau)^{-1}, \quad \mbox{for } 0\leq t-\tau \leq T.$$
	\end{enumerate}
	
\end{definition}

Conditions \ref{P_1} to \ref{P_3} ensure the existence of this family, as proved in \cite[Theorem 1]{sobol}. We briefly mention the ideas behind the construction of such family, since it depends on two auxiliary families of linear operators $\varphi_{\varepsilon}(t,\tau) \in \mathcal{L}(X)$ and $\Phi_{\varepsilon}(t,\tau) \in \mathcal{L}(X)$ that will be necessary in the sequel.

Suppose that $U_{\varepsilon}(t,\tau)\in \mathcal{L}(X)$ is a family satisfying the homogeneous differential equation given in \eqref{homogeneous_diff_eq}, that is, ${\partial_t} U_{\varepsilon}(t,\tau) = -A_{\varepsilon}(t) U_{\varepsilon}(t,\tau)$. Also, assume that there exists another family $\Phi_{\varepsilon} (t,\tau) \in \mathcal{L}(X)$ such that $U_{\varepsilon}(t,\tau)$ is obtained trough the integral equation
\begin{equation}\label{6.3}
U_{\varepsilon}(t,\tau) = T_{ -A_{\varepsilon}(\tau)  }(t-\tau) + \int_{\tau}^{t} T_{-A_{\varepsilon}(s)   } (t-s) \Phi_{\varepsilon} (s,\tau) ds.
\end{equation}

Differentiating in $t$, adding $A_{\varepsilon}(t) U_{\varepsilon}(t,\tau)$ on both sides and taking into account that ${\partial_t} U_{\varepsilon}(t,\tau) =-A_{\varepsilon}(t) U_{\varepsilon}(t,\tau) $, we deduce
\[
0  = \Phi_{\varepsilon}(t,\tau) - [A_{\varepsilon}(\tau)- A_{\varepsilon}(t)] T_{-A_{\varepsilon}(\tau) } (t-\tau) - \int_{\tau}^{t} [A_{\varepsilon}(s) - A_{\varepsilon}(t)] T_{-A_{\varepsilon}(s)} (t-s) \Phi_{\varepsilon} (s,\tau) ds.
\]

If we set 
\begin{equation}\label{varphi_1}
\varphi_{\varepsilon}(t,\tau) = [A_{\varepsilon}(\tau)- A_{\varepsilon}(t)] T_{-A_{\varepsilon}(\tau) } (t-\tau),
\end{equation}
then $\Phi_{\varepsilon}(t,\tau)$ would have to satisfy
\begin{equation}\label{6.6}
\Phi_{\varepsilon}(t,\tau) = \varphi_{\varepsilon}(t,\tau) + \int_{\tau}^{t} \varphi_{\varepsilon} (t,s) \Phi_{\varepsilon} (s,\tau) ds
\end{equation}
and it would be a fixed point of the map $S_{\varepsilon}(\Psi)(t) =\varphi_{\varepsilon}(t,\tau) + \int_{\tau}^{t} \varphi_{\varepsilon} (t,s) \Psi (s) ds$.

If we had a family $\Phi_{\varepsilon} (t,\tau)$ satisfying \eqref{6.6},
then we could proceed in the reverse way to obtain $U_{\varepsilon}(t,\tau)$. This is the technique employed to construct the linear process in the parabolic case \cite{sobol,Tanabe_1959} and the description of $U_{\varepsilon}(t,\tau)$ relies on this auxiliary family $\Phi_{\varepsilon}(t,\tau)$. The next proposition is proved in \cite[Section 5.6]{pazy} and \cite{sobol}. It ensures existence of  $\Phi_{\varepsilon}(t,\tau)$ and $U_{\varepsilon}(t,\tau)$ under the conditions required previously.

\begin{proposition}\label{ex_proc}
	For a fixed $\varepsilon \in [0,\varepsilon_0)$, assume that $\{ A_{\varepsilon}(t), \ t\in \mathbb{R}\}$ satisfies \ref{P_1}, \ref{P_2} and \ref{P_3}. Let $\delta \in (0,1]$ be the constant of H\"{o}lder continuity and $\{ \varphi_{\varepsilon} (t,\tau) \in \mathcal{L} (X), t \geq\tau \}$  the family given by \eqref{varphi_1}, then:
	
	\begin{enumerate}
		\item $\{  (t,\tau) \in \mathbb{R}^{2}; t > \tau  \} \ni (t,\tau) \mapsto
		\varphi_{\varepsilon}(t,\tau)
		\in \mathcal{L} (X)$ is continuous in the uniform	topology and
		\begin{equation*}%\label{vaphi_1_e}
		\| \varphi_{\varepsilon} (t,\tau) \|_{\mathcal{L}(X)} \leq C (t-\tau)^{\delta-1}, \quad \mbox{ for all } t > \tau, \tau \in \mathbb{R}.
		\end{equation*}
		
		\item 	There exists a unique family $\{ \Phi_{\varepsilon} (t,\tau) \in \mathcal{L} (X), t\geq \tau \}$ that satisfies \eqref{6.6} and this family is continuous in terms of the parameters $(t,\tau)$, that is, $\{  (t,\tau) \in \mathbb{R}^{2}; t> \tau  \} \ni (t,\tau) \mapsto
		\Phi_{\varepsilon}(t,\tau)
		\in \mathcal{L} (X)$ is continuous and for each $T>0$, there exists $C=C(T)>0$ such that 
		\begin{equation*}
		\| \Phi_{\varepsilon} (t,\tau) \|_{\mathcal{L} (X)} \leq C (t-\tau)^{\delta  -1}, \quad \mbox{ for all } 0< t-\tau \leq T.
		\end{equation*}
	\end{enumerate}

	Furthermore, the family of linear operators $\{U_{\varepsilon}(t,\tau) \in \mathcal{L}(X), t\geq \tau\}$ given by
	\[
	U_{\varepsilon}(t,\tau) = T_{-A_{\varepsilon}(\tau)} (t-\tau) + \int_{\tau}^{t} T_{-A_{\varepsilon}(s)} (t-s) \Phi_{\varepsilon}(s,\tau) ds
	\]
	is a \emph{linear process associated to $\{A_{\varepsilon}(t), \ t\in \R\}$} and satisfies the conditions in Definition \ref{evol_process}.
\end{proposition}

The fact that $U_{\varepsilon}(t,\tau)$ given by \eqref{6.3} satisfies all the conditions in Definition \ref{evol_process} can be found in the work of Sobolevski\u{\i} in \cite{sobol} or in the works of Tanabe \cite{Tanabe_1959, Tanabe_1960, Tanabe_1960_2}. Those four families of linear operators,  $T_{-A_{\varepsilon}(\tau)}(t-\tau)$, $U_{\varepsilon}(t,\tau)$, $\varphi_{\varepsilon} (t,\tau)$ and $\Phi_{\varepsilon}(t,\tau)$,  are essential to describe the  dynamics of the system associated to \eqref{Pe}. We are then able to enunciate one of the main results on this paper, the Toter-Kato type result on the convergence of the linear process as $\varepsilon \to 0^{+}$, whose proof is postponed to Section \ref{S_Estimates_rates_convergence}.

\begin{theorem}\label{T_Conv_Proc}
	Assume that conditions \ref{P_1} to \ref{P_5} hold, and let $\beta \in (0,1]$ be the constant in the resolvent estimate \eqref{res_est_X,Y}. For any $\theta \in (0,1)$, there exist constants $K,C>0$, independent of $\varepsilon \in [0,\varepsilon_0),$  such that
	\begin{align*}
	\| U_{\varepsilon}(t,\tau) - U_0(t,\tau) \|_{\mathcal{L}(X)} &\leq C(t-\tau)^{-\theta} e^{K(t-\tau)} \ell (\theta, \varepsilon) ,
	%\label{Conv_Proc_X}
	\\
	\| U_{\varepsilon}(t,\tau) - U_0(t,\tau) \|_{\mathcal{L}(X,Y)} &\leq C(t-\tau)^{-1+\beta (1-\theta)} e^{K(t-\tau)} \ell (\theta, \varepsilon), 
	%\label{Conv_Proc_X_Y}
	\end{align*}
	for all $\tau \in \mathbb{R}$ and $t>\tau$, where $\ell (\theta, \varepsilon) = \max \{ [\eta(\varepsilon)]^{\theta}, [\xi(\varepsilon)]^{\theta} \}$. In particular, $\ell (\theta, \varepsilon) \stackrel{\varepsilon \to 0^{+}}{\longrightarrow} 0.$
\end{theorem}

\subsection{Rate of convergence for the solution of the semilinear problem}

In order to obtain existence of global solution and convergence of them as $\varepsilon \to 0^{+}$, we need to require some properties on the family of nonlinearities $F_{\varepsilon}: \mathbb{R} \times Y \to X$. Assume that

\vspace{0.3cm}

\begin{enumerate}[label={\textbf{(NL.\arabic*)}},ref=(NL.\arabic*)]
	
	\item \label{(NL_1)}  Each $F_{\varepsilon} = F_{\varepsilon}(t,u)$ is H\"{o}lder continuous in $t$, globally Lipschitz in $u$ and bounded. Moreover, the constants $L>0$ of Lipschitz and $M>0$ of boundedness for $F_{\varepsilon}$ can be chosen uniformly in $\varepsilon$, that is
	\begin{align*}
	\| F_{\varepsilon} (t,u)\|_X 
	& \leq M, \quad \mbox{ for all } (t,u) \in \mathbb{R} \times Y, \ \varepsilon \in [0,\varepsilon_0), \\
	\| F_{\varepsilon} (t,u) - F_{\varepsilon} (t,v)\|_X 
	&\leq L \|u-v\|_Y, \quad \mbox{ for all }\varepsilon \in [0,\varepsilon_0), \ t\in \mathbb{R}, \ u,v \in Y.
	\end{align*}

	\item \label{(NL_2)} There exists a continuous  function $\gamma: [0,\varepsilon_0) \to \mathbb{R}^{+}$ with $\gamma(0)=0$ such that 
	\begin{equation*}
	\sup_{t\in \mathbb{R}} \sup_{u\in Y} \|F_{\varepsilon}(t,u) - F_0(t,u)\|_X \leq \gamma (\varepsilon).
	\end{equation*}

\end{enumerate}

\begin{remark}\label{R_cut-off}
		Conditions required in \ref{(NL_1)} are very restrictive and usually not found in practice. However, in many situations (like the application considered in Section \ref{S:Case_I_application_I}) we are able to prove that dynamics of the Problems \eqref{Pe} eventually enters a bounded subset of $Y$, uniformly in $\varepsilon$. If that is the case, we can proceed with a cut-off for the nonlinearities outside this bounded set so that the new family $F_{\varepsilon}$ obtained after the cut-off satisfies the assumptions required in \ref{(NL_1)}. The parabolic problems with this new nonlinearity will differ out-side the bounded set, but remains the same inside it, where all the solutions eventually go. Therefore, by restricting our attention to this uniform bounded absorbing set, we can assume that $F_{\varepsilon}$ have those desired properties inside it.
\end{remark}

It follows from \cite[Theorem 7]{sobol} that Problem \eqref{Pe} is locally well-posed, that is, there exists 
$$u_{\varepsilon}: [\tau, \tau + T(\varepsilon, \tau, u_{}^{\tau}) ) \to Y \mbox{ given by } u_{\varepsilon}(t) = U_{\varepsilon}(t,\tau) u_{}^{\tau} + \int_{\tau}^{t} U_{\varepsilon}(t,s) F_{\varepsilon}(s, u_{\varepsilon}(s)) ds,$$
solution of \eqref{Pe}, where $T(\varepsilon, \tau, u^{\tau})>0$ is the maximal interval of definition of $u_{\varepsilon}(t)$, and it depends on the initial condition and on $\varepsilon$. We denote the solution by $u_{\varepsilon}(t,\tau,u^{\tau})$ if we wish to emphasize the initial condition. 

In Section \ref{S_Estimates_rates_convergence}, Lemma \ref{L_est_sol}, we shall prove that the boundedness required for $F_{\varepsilon}$, implies that $\| u_{\varepsilon}(t)\|_Y$ remains bounded for $t$ in any interval of the form $[\tau, \tau+T]$. Therefore, the solution is globally defined in time and originates a  nonlinear process $\{S_{\varepsilon}(t,\tau): Y \to Y, \ t\geq \tau, \ \tau \in \mathbb{R} \}$ given by
\begin{equation*}\label{NL_process}
S_{\varepsilon}(t,\tau)u^{\tau} = u_{\varepsilon}(t,\tau,u^{\tau}).
\end{equation*}

	We now present the result on convergence of the solution as $\varepsilon \to 0^{+}$. Its proof is postponed to Section \ref{S_Estimates_rates_convergence}. 

\begin{theorem}\label{T_Conv_solution}
	Assume that conditions \ref{P_1} to \ref{P_5} hold, as well as \ref{(NL_1)} and \ref{(NL_2)}. Let $\beta\in (0,1]$ be the constant in the resolvent estimate \eqref{res_est_X,Y}. For any  $\theta \in (0,1)$, there exists constants $C,K>0$, independent of $\varepsilon, t , \tau$ such that, for any $\varepsilon\in [0,\varepsilon_0),$ $t>\tau$, $\tau\in \mathbb{R}$ and $u^{\tau }\in Y$, we have
	\begin{equation*}\label{Conv_solution}
	\| u_{\varepsilon}(t,\tau, u^{\tau}) - u_0 (t,\tau,u^{\tau}) \|_Y \leq C (t-\tau)^{-1+\beta(1-\theta)} e^{K(t-\tau)} \left[1+ \|u^{\tau}\|_Y\right]\rho(\theta, \varepsilon), 
	\end{equation*}
	where 
	\[
	\rho(\theta,\varepsilon) = \max\{ [\eta(\varepsilon)]^{\theta}, [\xi(\varepsilon)]^{\theta}, \gamma(\varepsilon)  \}.
	\]
	
	\noindent In particular, the convergence of the solution is uniform for $t$ in any interval of the form $[\tau+m,\tau+M]$, $0<m<M$, and $u^{\tau} \in B \subset Y$, $B$ bounded.

\end{theorem}

As an immediate consequence of the previous theorem, we have the following result.

\begin{corollary}\label{C_NL_process_convergence}
	Assume conditions \ref{P_1} to \ref{P_5}, \ref{(NL_1)} and \ref{(NL_2)} hold. Let $S_{\varepsilon}(t,\tau):Y \to Y$ be the nonlinear process obtained from the solution of \eqref{Pe}. For any compact set $I \subset (0,\infty)$ and any bounded set $B \subset Y$, we have
	\begin{equation*}
	\sup_{t\in I} \ \sup_{\tau \in \mathbb{R}} \  \sup_{u^{\tau} \in B} \| S_{\varepsilon}(t+\tau, \tau) u^{\tau}-S_0(t+\tau, \tau) u^{\tau}\|_Y \stackrel{\varepsilon \to 0}{\longrightarrow} 0. 
	\end{equation*}

\end{corollary}

\vspace{0.1cm}

\section{Estimates and rates of convergence}\label{S_Estimates_rates_convergence}

	This section is dedicated to obtain estimates and rate of convergence for a series of linear operators, culminating with the proof of the Trotter-Kato type Approximation result for the linear process $U_{\varepsilon}(t,\tau)$ to $U_0(t,\tau)$ and the convergence of the solution for the semilinear problem. We shall assume during this entire section that Conditions \ref{P_1} to \ref{P_5} and \ref{(NL_1)} to \ref{(NL_2)} hold.

\subsection{Resolvent convergence in $\mathcal{L}(X)$ and $\mathcal{L}(X,Y)$}

	We first estimate convergence of the resolvent of  $A_{\varepsilon}(t)$ to the resolvent of $A_0(t)$ in terms of $\varepsilon$. As an immediate consequence of  \ref{P_1} and \ref{P_2}, we obtain  existence of a constant $C>0$ (uniform in $t\in \mathbb{R}$ and $\varepsilon\in [0,\varepsilon_0)$) such that, for any $\lambda \in \Sigma_{\varphi} \cup\{0\}$, $\varepsilon\in [0,\varepsilon_0)$ and $t\in \mathbb{R}$, 
	\begin{equation}\label{est_A_res}
	\| A_{\varepsilon}(t) (\lambda + A_{\varepsilon}(t))^{-1} \|_{\mathcal{L}(X)} \leq C
	\
	\mbox{ and } 
	\
	\| A_{\varepsilon}(t) (\lambda + A_{\varepsilon}(t))^{-1} \|_{\mathcal{L}(Y)} \leq C.
	\end{equation}

	Since $Y \hookrightarrow X$, there exists a constant $C>0$ such that, for each $u\in Y$,
	\begin{equation}\label{cons_emb}
	\| u \|_X \leq C \| u \|_Y \quad \mbox{ and } \quad \| I\|_{\mathcal{L}(Y,X)} \leq C.
	\end{equation}

	Consequently, 
	\begin{align*}
	\| A_{\varepsilon}(t)^{-1} u - A_{0}(t)^{-1} u\|_X \leq C	\| A_{\varepsilon}(t)^{-1} u - A_{0}(t)^{-1} u\|_Y \leq C \eta (\varepsilon) \|u\|_X 
	\end{align*}
	and
	\begin{align*}
	\| A_{\varepsilon}(t)^{-1} u - A_{0}(t)^{-1} u\|_Y \leq  \eta (\varepsilon) \|u\|_X \leq C \eta (\varepsilon) \|u\|_Y.
	\end{align*}
	
	Hence, we have the following estimates 
	\begin{equation}\label{inv_ope_conv_Lx_Ly}
	 \| A_{\varepsilon}(t)^{-1}  - A_{0}(t)^{-1} \|_{\mathcal{L}(X)} \leq C \eta (\varepsilon)
	 \quad \mbox{ and } \quad 
	 \| A_{\varepsilon}(t)^{-1}  - A_{0}(t)^{-1} \|_{\mathcal{L}(Y)} \leq C \eta (\varepsilon)
	\end{equation}

	From the resolvent equality and simple algebra we can prove that the following equalities hold, for all  $\lambda \in \Sigma_{\varphi}\cup \{0\}$ and $t\in \mathbb{R}$,
	\begin{align}
	\scalebox{0.95}{
	$
	(\lambda + A_{\varepsilon}(t))^{-1} - (\lambda+A_{0}(t))^{-1} $} &
	\scalebox{0.95}{$ = A_{\varepsilon}(t) (\lambda+A_{\varepsilon}(t))^{-1}  [ A_{\varepsilon}(t)^{-1} - A_{0}(t)^{-1}]A_{0}(t) (\lambda+A_{0}(t))^{-1},
	$}
	\label{res_dif} \\
	\scalebox{0.95}{
	$A_{\varepsilon}(t) (\lambda+A_{\varepsilon}(t))^{-1}-A_{0}(t)(\lambda+A_{0}(t))^{-1}$} & \scalebox{0.95}{$= - \lambda( \lambda+A_{\varepsilon}(t) )^{-1}A_{\varepsilon}(t) [ A_{\varepsilon}(t)^{-1} - A_0(t)^{-1}  ] A_{0}(t)  (\lambda+A_{0}(t))^{-1}.
	$}
	\label{res_dif_2} 
	\end{align}

\

	Expression \eqref{res_dif} implicates that resolvent convergence inside the sector $\Sigma_{\varphi} \cup \{0\}$  follows from the convergence of $ A_{\varepsilon}(t)^{-1} - A_{0}(t)^{-1}$, as $\varepsilon \to 0^{+}$, requested in \ref{P_5}, as stated in next proposition.

\begin{proposition}\label{P_resolvent_convergence}
	There exists a constant $C>0$, independent of $\varepsilon \in [0,\varepsilon_0)$ or $t\in \mathbb{R}$, such that, for all $\lambda \in \Sigma_{\varphi} \cup\{0\}$ and $t\in \mathbb{R}$, 
	\begin{align*}
	  \| (\lambda+A_{\varepsilon}(t))^{-1} - (\lambda+A_{0}(t))^{-1} \|_{\mathcal{L} (X,Y)} & \leq C \eta (\varepsilon).
	%\label{res_conv_X_Y}
	\end{align*}
	
	\begin{proof}
		From the uniform estimate obtained in \eqref{est_A_res} for  $A_{\varepsilon}(t) (\lambda +A_{\varepsilon}(t))^{-1}$ in  $\mathcal{L}(X)$ and $\mathcal{L}(Y)$ and Equality \eqref{res_dif}, we deduce
		\begin{align*}
		& \|  (\lambda+A_{\varepsilon}(t))^{-1} - (\lambda+A_{0}(t))^{-1}      \|_{\mathcal{L} (X,Y)} \\
		& \quad \leq 
		\| A_{\varepsilon}(t) (\lambda+A_{\varepsilon}(t))^{-1}  \|_{\mathcal{L}(Y)} 
		\| A_{\varepsilon}(t)^{-1}   - A_{0}(t)^{-1} \|_{\mathcal{L}(X,Y)} 
		\| A_{0}(t) (\lambda+A_{0}(t))^{-1} \|_{\mathcal{L}(X)}\\
		& \quad \leq
		C\eta (\varepsilon).
		\end{align*}		
	\end{proof}
\end{proposition}

Another estimate on the resolvent in terms of $\varepsilon$ that will be useful in the sequel is presented next.

\begin{lemma}\label{L_est_A_res}
	There exists a constant $C>0$, independent of $\varepsilon \in (0,\varepsilon_0]$ and $t\in \mathbb{R}$, such that,  , 
	\begin{equation*}
	\| A_{\varepsilon} (t) (\lambda + A_{\varepsilon}(t))^{-1} - A_{0}(t) (\lambda+A_0(t))^{-1} \|_{\mathcal{L}(X)} \leq C |\lambda |\eta( \varepsilon), \quad \mbox{ for any } \lambda \in \Sigma_{\varphi} \cup\{0\}.
	\end{equation*}
	
			\begin{proof}
				It follows directly from the estimates \eqref{est_A_res} for  $A_{\varepsilon}(t) (\lambda +A_{\varepsilon}(t))^{-1}$ in  $\mathcal{L}(X)$ and $\mathcal{L}(Y)$, from \eqref{cons_emb} and from \eqref{res_dif_2} that
				\begin{align*}
				& \|  A_{\varepsilon}(t) (\lambda+A_{\varepsilon}(t))^{-1} - A_{0}(t) (\lambda+A_{0}(t))^{-1}      \|_{\mathcal{L} (X)} \\
				& \quad \leq 
				\|\lambda A_{\varepsilon}(t) (\lambda+A_{\varepsilon}(t))^{-1}  \|_{\mathcal{L}(Y,X)} 
				\| A_{\varepsilon}(t)^{-1}   - A_{0}(t)^{-1} \|_{\mathcal{L}(X,Y)} 
				\| A_{0}(t) (\lambda+A_{0}(t))^{-1} \|_{\mathcal{L}(X)}\\
				& \quad \leq 
				|\lambda | \|I\|_{\mathcal{L}(Y,X)} \|A_{\varepsilon}(t) (\lambda+A_{\varepsilon}(t))^{-1}  \|_{\mathcal{L}(Y)} 
				\| A_{\varepsilon}(t)^{-1}   - A_{0}(t)^{-1} \|_{\mathcal{L}(X,Y)} 
				\| A_{0}(t) (\lambda+A_{0}(t))^{-1} \|_{\mathcal{L}(X)}\\
				& \quad \leq
				C|\lambda| \eta (\varepsilon).
				\end{align*}
			\end{proof}
\end{lemma}

	Lastly on the linear operator and its resolvent, we provide an estimate for a situation where we vary both $\varepsilon$ and time $t\in \mathbb{R}$ simultaneously.

	\begin{lemma}\label{L_est_A_e_t}
		Let $\delta \in (0,1]$ be the H\"{o}lder continuity constant in \ref{P_3}.
		For any $\theta \in [0,1]$, there exists a constant $C>0$ such that, for all  $t,\tau \in \mathbb{R}$ and $\varepsilon\in [0,\varepsilon_0),$ 
		\begin{equation}\label{est_A_e_t}
		\| A_{\varepsilon}(t) A_{\varepsilon}(\tau)^{-1} - A_{0}(t) A_0(\tau)^{-1}\|_{\mathcal{L}(X)} 
		\leq C |t-\tau|^{\delta (1-\theta)} [\xi (\varepsilon)]^{\theta}.
		\end{equation}
		
		\begin{proof}
			From \ref{P_3}, we deduce
			\begin{align}
			\| A_{\varepsilon}(t) A_{\varepsilon}(\tau)^{-1} - A_0(t) A_0(\tau)^{-1} \|_{\mathcal{L}(X)} 
			& = 	\| [ A_{\varepsilon} (\tau) - A_{\varepsilon}(t)] A_{\varepsilon}(\tau)^{-1} - [A_{0}(\tau) -A_0(t)] A_0(\tau)^{-1} \|_{\mathcal{L}(X)} \nonumber \\ 
			& \leq C|t-\tau|^{\delta}. \label{est_1}
			\end{align}
			
			Now, interpolating \eqref{est_1} and the estimate in \ref{P_4} with an exponent $\theta \in [0,1]$, we obtain \eqref{est_A_e_t}.
		\end{proof}
	\end{lemma}

\

\subsection{Convergence and estimates for the semigroups}

	Since each operator $A_{\varepsilon}(\tau)$ is sectorial (with an uniform sector and uniform resolvent estimates in terms of $\varepsilon$ and $\tau$), classical theory on  semigroups implies that $-A_{\varepsilon}(\tau)$ generates an analytic semigroup, which we denote by $T_{-A_{\varepsilon}(\tau)}(\cdot)$.
	
	If  $\Gamma$ is the contour of $\Sigma_{\varphi} \subset \rho ( -A_{\varepsilon}(\tau) )$,  that is, $\Gamma = \{re^{- i \varphi}: r>0\} \cup \{re^{ i \varphi}: r>0\} $ and it is oriented with increasing imaginary part, then we have the following expressions 
	\begin{align}
	T_{-A_{\varepsilon}(\tau)} (t) & = \frac{1}{2\pi i} \int_{\Gamma}  e^{\lambda t} (\lambda +A_{\varepsilon}(\tau))^{-1} d\lambda,  \label{sem_exp}\\ 
	A_{\varepsilon} (\tau) T_{-A_{\varepsilon}(\tau)} (t) & = \frac{1}{2\pi i} \int_{\Gamma}  e^{\lambda t} A_{\varepsilon} (\tau) (\lambda +A_{\varepsilon}(\tau))^{-1} d\lambda,  \label{A_sem_exp}
	\end{align}
	that can be found in \cite[Section 2.5]{pazy}. A direct application of estimates \eqref{res_est_X} and \eqref{est_A_res} in Expressions \eqref{sem_exp} and \eqref{A_sem_exp} implies, for any $\tau \in \mathbb{R}$ and $\varepsilon \in [0,\varepsilon_0)$, 
	\begin{align}
	& \|T_{-A_{\varepsilon}(\tau)}(t) \|_{\mathcal{L}(X)} \leq C, \quad \mbox{ for all }t\geq 0, \label{Teorema_2.5.2_(A)} \\
	& \| A_{\varepsilon} (\tau) T_{-A_{\varepsilon}(\tau)} (t) \|_{\mathcal{L}(X)} \leq Ct^{-1}, \quad \mbox{ for all }t> 0 . \label{Pazy_2.5.6}
	\end{align}
	
	Uniformity of those estimates with respect to $\varepsilon$ and $\tau$ follows from \ref{P_1} and \ref{P_2}. We can also obtain an estimate for this semigroup in $\mathcal{L}(X,Y)$, as stated next.

	\begin{lemma}\label{L_Sem_Est_X_Y}
		
			Let $\beta \in (0,1]$ be the constant in \ref{P_2}. There exists a constant $C>0$ independent of $\varepsilon \in [0,\varepsilon_0)$ and  $\tau \in \mathbb{R}$, such that, for all $t>0$,
		\begin{equation*}
		\| T_{-A_{\varepsilon}(\tau)} (t) \|_{\mathcal{L}(X,Y)} \leq C t^{\beta-1}.
		\end{equation*}
		\begin{proof}
			Using estimate \eqref{res_est_X,Y}, we obtain
			\begin{align*}
			\| T_{-A_{\varepsilon}(\tau)} (t) \|_{\mathcal{L}(X,Y)} 
			& \leq \frac{1}{2\pi}  \int_{\Gamma} |e^{\lambda t} | \| (\lambda+A_{\varepsilon}(\tau))^{-1}\|_{\mathcal{L}(X,Y)} |d\lambda| 	
			\leq C  \int_{0}^{\infty} e^{r [\cos \varphi] t } \frac{C}{1+r^{\beta}} dr \\
			& \leq 
			Ct^{\beta -1} \int_{0}^{\infty} e^{[\cos \varphi] u } \frac{1}{t^{\beta}+ u^{\beta}} du 
			=
			C(\varphi, \beta) t^{\beta -1},
			\end{align*}
			where constant $C$ depends on the angle $\varphi$ and on $\beta$, but it is independent of $\varepsilon,\tau$ and $t$. 
		\end{proof}
	\end{lemma}

	We establish next a convergence of the linear semigroups relative to $\varepsilon$.

	\begin{lemma}
	Let $\beta \in (0,1]$ be the constant in \ref{P_2}. For any $\theta \in [0,1]$, there exists a  constant $C>0$ independent of $\varepsilon \in [0,\varepsilon_0)$ and  $\tau \in \mathbb{R}$, such that, for all $t>0$, 
	\begin{align}
	&	\| T_{-A_{\varepsilon}(\tau)} (t) - T_{-A_0(\tau)}(t)  \|_{\mathcal{L}(X)} \leq C t^{-\theta} [\eta(\varepsilon)]^{\theta} \label{conv_sem_e_X}, \\
	&\| T_{-A_{\varepsilon}(\tau)} (t) - T_{-A_0(\tau)}(t)  \|_{\mathcal{L}(X,Y)} \leq C t^{-1+\beta(1-\theta)} [\eta(\varepsilon)]^{\theta}. \label{conv_sem_e_X,Y} 
	\end{align}

		\begin{proof}
		It follows from \eqref{Teorema_2.5.2_(A)} that 
		\begin{equation}\label{sem_dif}
		\| T_{-A_{\varepsilon}(\tau)}(t) - T_{-A_0(\tau) }(t)  \|_{\mathcal{L} (X) } \leq C.
		\end{equation}
		
		Consider the curve $\Gamma$ parametrized as  $\Gamma = \Gamma_1 \vee \Gamma_2^{-}$ where 
		\begin{align*}
		& \Gamma_1 : = \{ \lambda \in \mathbb{C}: \ \lambda =  re^{i \varphi};\  r\in [0,\infty)  \}, 
		\quad  \Gamma_2 : = \{ \lambda \in \mathbb{C}: \ \lambda = re^{-i \varphi}; \  r\in [0,\infty)  \},
		\end{align*}
		and $\Gamma_2^{-}$ stands for the reverse path. Using the symmetry of curves $\Gamma_1$ and $\Gamma_2$ and estimate \eqref{inv_ope_conv_Lx_Ly},  we obtain
		\begin{align}
		\| T_{-A_{\varepsilon}(\tau)}(t) - T_{-A_0(\tau) }(t)   \|_{\mathcal{L}(X)} 
		& \leq \frac{1}{\pi} \int_{\Gamma_1} |e^{\lambda t}| \|   (\lambda+A_{\varepsilon}(\tau))^{-1} - (\lambda+A_0(\tau))^{-1} \|_{\mathcal{L}(X)} | d\lambda| 
		\nonumber\\
		& \leq C\eta(\varepsilon) \int_{0}^{\infty}   e^{r [\cos \varphi] t} dr
		\nonumber\\
		&  \leq C({\varphi}) t^{-1}  \eta (\varepsilon)	\label{sem_dif_2}
		\end{align}
		
		 Interpolating \eqref{sem_dif} and \eqref{sem_dif_2} with exponents $\theta$ and $1-\theta$, for $\theta \in [0,1]$,  we obtain  \eqref{conv_sem_e_X}. In order to estimate \eqref{conv_sem_e_X,Y}, we first note from Lemma \ref{L_Sem_Est_X_Y} that 
		\begin{equation}\label{sem_dif_X_Y}
		\| T_{-A_{\varepsilon}(\tau)}(t) - T_{-A_0(\tau) }(t)  \|_{\mathcal{L} (X,Y) } \leq Ct^{\beta-1}.
		\end{equation}

		Using \ref{P_5} and the integral formulation for the semigroup \eqref{semigrupo}, we obtain
		\begin{align}
		\| T_{-A_{\varepsilon}(\tau)}(t) - T_{-A_0(\tau) }(t)   \|_{\mathcal{L}(X,Y)} 
		& \leq \frac{1}{\pi} \int_{\Gamma_1} |e^{\lambda t}| \|   (\lambda+A_{\varepsilon}(\tau))^{-1} - (\lambda+A_0(\tau))^{-1} \|_{\mathcal{L}(X,Y)} | d\lambda| 
		\nonumber\\
		&  \leq C_{\varphi} t^{-1}  \eta (\varepsilon)	\label{sem_dif_X_Y_2}
		\end{align}
		
		Interpolating \eqref{sem_dif_X_Y} and \eqref{sem_dif_X_Y_2} with exponents $1-\theta$ and $\theta$,  $\theta \in [0,1]$, we obtain the desired estimate \eqref{conv_sem_e_X,Y}.
		
	\end{proof}
	\end{lemma}

	We deduce in the sequel the last estimate on semigroup in terms of $\varepsilon$ necessary to our future analysis.

	\begin{lemma}\label{L_A_sem}
		For any $\theta \in [0,1]$, there exists $C>0$ such that, for all $\varepsilon \in [0,\varepsilon_0),$ $\tau \in \mathbb{R}$  and $t>0$,
		\begin{equation*}
		\| A_{\varepsilon}(\tau) T_{-A_{\varepsilon}(\tau)} (t) -A_{0} (\tau) T_{-A_0(\tau)}(t) \|_{\mathcal{L}(X)} \leq C t^{-1-\theta} [\eta (\varepsilon)]^{\theta}.
		\end{equation*}
		
	\begin{proof}
		Note that from \eqref{Pazy_2.5.6}, we obtain
		\begin{equation}\label{A_sem_est}
		\| A_{\varepsilon}(\tau) T_{-A_{\varepsilon}(\tau)} (t) -A_{0} (\tau) T_{-A_0(\tau)}(t) \|_{\mathcal{L}(X)} \leq C t^{-1}.
		\end{equation}
		
		On the other hand, using the integral formulation for the semigroup and Lemma \ref{L_est_A_res}, we deduce
		\begin{align}
		 \| A_{\varepsilon}(\tau) T_{-A_{\varepsilon}(\tau)} (t) -A_{0} (\tau) T_{-A_0(\tau)}(t) \|_{\mathcal{L}(X)} 
		  & \leq C \int_{\Gamma} |e^{\lambda t}| 
		\| 
			A_{\varepsilon}(\tau) (\lambda +A_{\varepsilon}(\tau) )^{-1} - A_0(\tau) (\lambda + A_0(\tau) )^{-1} 
		\|_{\mathcal{L}(X)} |d\lambda| \nonumber\\
		&   \leq C \eta(\varepsilon) \int_{0}^{\infty} e^{r [\cos \varphi] t} r dr
		  \leq C(\varphi) t^{-2} \eta(\varepsilon). \label{A_sem_est_2} 
		\end{align}
		
		Interpolating \eqref{A_sem_est} and \eqref{A_sem_est_2} with exponents $1-\theta$ and $\theta$, we obtain the desired estimate.
	\end{proof}
	\end{lemma}

\subsection{Convergence and estimates for the families $\varphi_{\varepsilon}(t,\tau)$ and $\Phi_{\varepsilon}(t,\tau)$}

In order to achieve our final goal of obtaining rate of convergence for the linear process associated to the family $\{A_{\varepsilon}(t), \ t\in \mathbb{R}\}$, we first need to establish rate of convergences for the auxiliary families $\varphi_{\varepsilon}(t,\tau)$ and $\Phi_{\varepsilon}(t,\tau)$. Recall that 
\[
\varphi_{\varepsilon}(t,\tau) = [A_{\varepsilon}(\tau) - A_{\varepsilon}(t)] T_{-A_{\varepsilon}(\tau)} (t-\tau),
\]
and it follows directly from \ref{P_3} and \eqref{Pazy_2.5.6} that 
\[
\| \varphi_{\varepsilon} (t,\tau) \|_{\mathcal{L}(X)}\leq C(t-\tau)^{\delta-1}, \quad \mbox{ for any }t>\tau.
\]

The rate of convergence required for the resolvent operators in Properties \ref{P_4} and \ref{P_5} are transfered to the families $\varphi_{\varepsilon}(t,\tau)$ as follows.

\begin{lemma}\label{L_est_varphi}
	Let $\theta \in [0,1]$ and $\delta$ be the constant of H\"{o}lder continuity in \ref{P_3}. There exists a constant $C>0$ such that, for any $\varepsilon \in [0,\varepsilon_0)$, $t> \tau $ and $\tau \in \mathbb{R}$, we have
	\begin{equation*}\label{est_varphi}
	\| \varphi_{\varepsilon}(t,\tau) - \varphi_0 (t,\tau)\|_{\mathcal{L}(X)} \leq C (t-\tau)^{-1+\delta (1-\theta)} \ell (\theta, \varepsilon),
	\end{equation*}
	where 
	\begin{equation}\label{l_theta_e}
	\ell(\theta, \varepsilon) = \max \{  [\eta (\varepsilon)]^{\theta} , [\xi (\varepsilon)]^{\theta}  \}. 
	\end{equation}

	\noindent In particular,  $\ell (\theta,\varepsilon) \stackrel{\varepsilon \to 0}{\longrightarrow} 0$.
	
	\begin{proof}
	Using the previous estimates and Expression \eqref{varphi_1} for the family $\varphi_{\varepsilon} (t,\tau)$, we obtain
	\begin{align*}
	& \| \varphi_{\varepsilon} (t,\tau) - \varphi_0(t,\tau) \|_{\mathcal{L}(X)} \\
	%
	%
	%Primeira desigualdade
	%
	%
	& \quad \leq 
	\|
	 [A_{\varepsilon}(\tau) -A_{\varepsilon}(t)] A_{\varepsilon}(\tau)^{-1} A_{\varepsilon}(\tau) T_{-A_{\varepsilon}(\tau)} (t-\tau) 
	 -
	 [A_0(\tau)-A_0(t)] A_0(\tau)^{-1}A_0(\tau)T_{-A_0(\tau)}(t-\tau)
	\|_{\mathcal{L}(X)} \\
	%	
	%
	%Segunda Desigualdade
	%
	%
	& \quad \leq 
	\|[A_{\varepsilon}(\tau)-A_{\varepsilon}(t)] A_{\varepsilon}(\tau)^{-1} \|_{\mathcal{L}(X)} 
	\|
	A_{\varepsilon}(\tau) T_{-A_{\varepsilon}(\tau)} (t-\tau) -A_0(\tau)T_{-A_0(\tau)}(t-\tau)
	\|_{\mathcal{L}(X)}\\
	& \qquad + 
	\|
		A_{\varepsilon}(t)A_{\varepsilon}(\tau)^{-1} - A_0(t)A_0(\tau)^{-1} 
	\|_{\mathcal{L}(X)}  
	\|
		A_0(\tau)T_{-A_0(\tau)}(t-\tau)
	\|_{\mathcal{L}(X)}\\
	%
	%
	%Terceira Desigualdade
	%
	%
	& \quad \leq C(t-\tau)^{-1+\delta-\theta} [\eta(\varepsilon)]^{\theta} + C(t-\tau)^{-1+\delta(1-\theta)} [ \xi(\varepsilon)]^{\theta} \\
	%
	%
	%Quarta desigualdade
	%
	%
	& \quad \leq C (t-\tau)^{-1+\delta (1-\theta)} \ell (\theta, \varepsilon),
	\end{align*}		
	since $-1+\delta(1-\theta) < -1+\delta -\theta<0.$ 
	\end{proof}
\end{lemma}

	As far as estimates for the family $\Phi_{\varepsilon} (t,\tau)$, we have the following result.
	
	\begin{lemma}\label{L_est_Phi}
	Let $\delta \in (0,1]$ be the constant of H\"{o}lder continuity in \ref{P_3}. There exist constants $C,K>0$ such that, for any $\varepsilon \in [0,\varepsilon_0)$, $\tau\in \mathbb{R}$ and $t>\tau$, 
	\begin{equation*}\label{est_Phi}
	\| \Phi_{\varepsilon}(t,\tau) \|_{\mathcal{L}(X)} \leq C(t-\tau)^{\delta-1} e^{K(t-\tau)}.
	\end{equation*}
	\begin{proof}
	From previous estimates, we obtain
	\begin{align*}
	\| \Phi_{\varepsilon}(t,\tau) \|_{\mathcal{L}(X)} & 
	\leq 	\| \varphi_{\varepsilon}(t,\tau) \|_{\mathcal{L}(X)} + \int_{\tau}^{t} 	\| \varphi_{\varepsilon}(t,s) \|_{\mathcal{L}(X)}	\| \Phi_{\varepsilon}(s,\tau) \|_{\mathcal{L}(X)}ds \\
	& \leq C(t-\tau)^{\delta-1} + \int_{\tau}^{t} C(t-s)^{\delta-1} 	\| \Phi_{\varepsilon}(s,\tau) \|_{\mathcal{L}(X)} ds.
	\end{align*}
	
	It follows from Gronwall's inequality \cite[p.190]{henry} that, for any $t>\tau$,
	\[
		\| \Phi_{\varepsilon}(t,\tau) \|_{\mathcal{L}(X)} \leq \frac{C}{\delta} (t-\tau)^{\delta-1} e^{K(t-\tau)},
	\]
	where $K>(2C\Gamma(\delta))^{\frac{1}{\delta}}$.
	\end{proof}
	\end{lemma}

	We obtain in the sequel a rate of convergence for the family $\Phi_{\varepsilon}$ as $\varepsilon$ converges to zero. Unlike previous results, our auxiliary $\theta$ that appears in the estimate needs to be in the open interval $(0,1)$ instead of $[0,1]$ in order to ensure convergence of integrals that feature in the estimates.

\begin{lemma}\label{L_Conv_Phi}
	Let $\theta \in (0,1)$ and $\delta$ be the constant of H\"{o}lder continuity in \ref{P_3}. There exists $C=C(\theta,\delta)>0$ and $K=K(\delta)>0$ such that, for any $\varepsilon \in [0,\varepsilon_0)$, $\tau \in \mathbb{R}$ and $t>\tau$, we have
	\begin{equation}\label{Conv_Phi}
	\| \Phi_{\varepsilon} (t,\tau) - \Phi_0(t,\tau) \|_{\mathcal{L}(X)} 
	\leq 
	C(t-\tau)^{-1+\delta(1-\theta)} e^{K(t-\tau)} \ell (\theta,\varepsilon),
	\end{equation}
	where $\ell (\theta, \varepsilon)$ is given in \eqref{l_theta_e}.

	\begin{proof}
		Using the estimates obtained earlier, we deduce
	\begin{align*}
		\| \Phi_{\varepsilon} (t,\tau) - \Phi_0(t,\tau) \|_{\mathcal{L}(X)} 
		%
		%
		%Primeira desigualdade
		%
		%
		& \leq \| \varphi_{\varepsilon} (t, \tau) - \varphi_0 (t,\tau) \|_{\mathcal{L}(X)} 
		+
		\int_{\tau}^{t} \| \varphi_{\varepsilon} (t,s) - \varphi_0 (t,s) \|_{\mathcal{L}(X)} \| \Phi_{\varepsilon}(s,\tau) \|_{\mathcal{L}(X)}
		ds\\
		& \quad  + \int_{\tau}^{t} \| \varphi_0 (t,s) \|_{\mathcal{L}(X)} \| \Phi_{\varepsilon}(s,\tau) - \Phi_0 (s,\tau) \|_{\mathcal{L}(X)} 
		ds\\
		%
		%
		%Segunda Desigualdade
		%
		%
		& \leq C(t-\tau)^{-1+\delta(1-\theta)}\ell(\theta, \varepsilon) 
		+
		\int_{\tau}^{t} C(t-s)^{-1+\delta(1-\theta)} \ell(\theta, \varepsilon)(s-\tau)^{\delta-1}
		ds\\
		& \quad  + \int_{\tau}^{t} C(t-s)^{\delta-1} \|\Phi_{\varepsilon}(s,\tau) - \Phi_0 (s,\tau) \|_{\mathcal{L}(X)}
		ds \\
		%
		%
		%Terceira desigualdade
		%
		%
		& \leq C(t-\tau)^{-1+\delta(1-\theta)}\ell(\theta, \varepsilon) 
		+
		 C(t-\tau)^{\delta(1-\theta)+\delta-1}\mathcal{B}(\delta(1-\theta), \delta) \ell(\theta, \varepsilon)
		\\
		& \quad  + \int_{\tau}^{t} C (t-s)^{\delta-1} \|\Phi_{\varepsilon}(s,\tau) - \Phi_0 (s,\tau) \|_{\mathcal{L}(X)} ,
	\end{align*}
	where $\mathcal{B}(\cdot, \cdot)$ is the Beta function. Taking $\psi (t) = \| \Phi_{\varepsilon} (t,\tau) - \Phi_0(t,\tau) \|_{\mathcal{L}(X)} $, we restate the above inequality as 
	\[
	\psi(t) \leq C \left[ (t-\tau)^{-1+\delta(1-\theta)} + (t-\tau)^{\delta(1-\theta)+\delta-1} \right]\ell(\theta, \varepsilon) + C\int_{\tau}^{t} (t-s)^{\delta-1}\psi(s)ds.
	\]
	
	Applying the generalized version of Gronwall's inequality \cite[p.190]{henry}, 
	\[
	\psi(t) \leq C(\delta, \theta) \left[(t-\tau)^{-1+\delta(1-\theta)}+(t-\tau)^{\delta(1-\theta)+\delta-1}\right] \ell(\theta,\varepsilon) e^{{K}(t-\tau)},
	\]
	for $K>(C\Gamma(\delta))^{\frac{1}{\delta}}.$ Moreover, if \scalebox{0.9}{$\delta(1-\theta) + \delta-1>0$}, then the growth provided by the term \scalebox{0.9}{$(t-\tau)^{\delta(1-\theta) + \delta-1}$} can be incorporated to the exponential term $e^{K(t-\tau)}$, correcting the constant if necessary. If \scalebox{0.9}{$\delta(1-\theta) + \delta-1<0$}, then \scalebox{0.9}{$(t-\tau)^{\delta(1-\theta) + \delta-1} \leq (t-\tau)^{\delta(1-\theta)-1}$} for $t-\tau$ near zero. In both cases, Inequality \eqref{Conv_Phi} follows from the above estimate.
	\end{proof}
\end{lemma}

\subsection{Convergence and estimates for the linear process $U_{\varepsilon}(t,\tau)$}

Before we prove Theorem \ref{T_Conv_Proc}, we obtain an estimate for the linear process that will be necessary.

\begin{lemma}
	Let $\beta \in (0,1]$ be the constant in \ref{P_2}. There exists $C,K>0$ such that, for any $\varepsilon \in [0,\varepsilon_0)$, $\tau \in \mathbb{R}$ and $t>\tau$,
	\begin{align*}
	\| U_{\varepsilon}(t,\tau) \|_{\mathcal{L}(X)} & \leq C  e^{K(t-\tau)}, \\
	\| U_{\varepsilon}(t,\tau) \|_{\mathcal{L}(X,Y)} & \leq C  (t-\tau)^{\beta-1} e^{K(t-\tau)}. 
	\end{align*}
	
\begin{proof}
From previous estimates and the expression for the linear process in \eqref{6.3}, we obtain
\begin{align*}
\|U_{\varepsilon}(t,\tau) \|_{\mathcal{L}(X)} 
&\leq \|T_{-A_{\varepsilon}(\tau)} (t-\tau) \|_{\mathcal{L}(X)} 
+ \int_{\tau}^{t} \|T_{-A_{\varepsilon}(s)} (t-s)\|_{\mathcal{L}(X)} \| \Phi_{\varepsilon} (s,\tau) \|_{\mathcal{L}(X) } ds \\
&\leq C+ \int_{\tau}^{t} C (s-\tau)^{\delta-1} e^{K(s-\tau)} ds \leq C+\frac{C}{\delta}(t-\tau)^{\delta}e^{K(t-\tau)} \\
&\leq C \max \{1,(t-\tau)^{\delta}\}e^{K(t-\tau)} \leq C e^{K(t-\tau)},
\end{align*}	
where we incorporated $(t-\tau)^{\delta}$ into the growth presented by $e^{K(t-\tau)}$, making adjustments in the constant $C$, if necessary. Similarly, we have
\begin{align*}
\|U_{\varepsilon}(t,\tau) \|_{\mathcal{L}(X,Y)} 
&
\leq C(t-\tau)^{\beta-1}+
 \int_{\tau}^{t} C(t-s)^{\beta-1} (s-\tau)^{\delta-1} e^{K(s-\tau)}ds  \\
&\leq C(t-\tau)^{\beta-1}+\frac{C}{\beta+\delta}(t-\tau)^{\beta+\delta-1} \mathcal{B}(\beta,\delta) e^{K(t-\tau)} \\
&\leq C(t-\tau)^{\beta-1}\left[1+(t-\tau)^{\delta} e^{K(t-\tau)} \right]\\
& \leq C (t-\tau)^{\beta-1} e^{K(t-\tau)},
\end{align*}
also adjusting the constant $C$, if necessary.
\end{proof}
\end{lemma}

We are now able to prove Theorem \ref{T_Conv_Proc}.

\

\noindent\textbf{Proof of Theorem \ref{T_Conv_Proc}: }
%\begin{proof}
		We first obtain the estimate in $\mathcal{L}(X)$ using Expression \eqref{6.3} for the linear process and the estimates established previously.
		\begin{align*}
		& \| U_{\varepsilon}(t,\tau) - U_0(t,\tau) \|_{\mathcal{L}(X)} \\
		%
		%
		%Primeira Desigualdade
		%
		%
		& \quad \leq  \| T_{-A_{\varepsilon} (\tau) } (t-\tau) - T_{-A_0(\tau)} (t-\tau)  \|_{\mathcal{L}(X)}
		+
		\int_{\tau}^{t} 
		\| 
		T_{-A_{\varepsilon}(s)}(t-s) 	\|_{\mathcal{L}(X)}
		\|
		[ \Phi_{\varepsilon} (s,\tau)-\Phi_0(s,\tau)]
		\|_{\mathcal{L}(X)}
		ds \\
		& \qquad +
		\int_{\tau}^{t}
		\| 
		[T_{-A_{\varepsilon}(s)} (t-s) - T_{A_{0} (s)} (t-s)] 
		\|_{\mathcal{L}(X)}
		\|
		\Phi_0(s,\tau)
		 \|_{\mathcal{L}(X)}
		 ds\\
		 %
		 %
		 %Segunda Desigualdade
		 %
		 %
		 & \quad \leq C(t-\tau)^{-\theta} [\eta(\varepsilon)]^{\theta}
		 +
		 \int_{\tau}^{t} 
		 C
		 (s-\tau)^{-1+\delta(1-\theta)} e^{K(s-\tau)} \ell (\theta, \varepsilon)
		 ds 
		  +
		 \int_{\tau}^{t}
		C(t-s)^{-\theta} [\eta(\varepsilon)]^{\theta}
		(s-\tau)^{\delta-1}e^{K(s-\tau)}
		 ds\\
		 %
		 %
		 %Terceira Desigualdade
		 %
		 %
		 & \quad \leq C(t-\tau)^{-\theta} [\eta(\varepsilon)]^{\theta}
		 +
		 \frac{C}{\delta(1-\theta)} e^{K(t-\tau)} \ell (\theta, \varepsilon)
		 (t-\tau)^{\delta(1-\theta)} 
		 +
		 C [\eta(\varepsilon)]^{\theta}
		 (t-\tau)^{\delta-\theta} \mathcal{B}(1-\theta,\delta)e^{K(t-\tau)}\\
		 %
		 %
		 %Quarta Desigualdade
		 %
		 %
		 & \quad \leq C(\theta, \delta) (t-\tau)^{-\theta} e^{K(t-\tau)} \ell(\theta, \varepsilon).
		\end{align*}
		
		Proceeding similarly to the estimate in $\mathcal{L}(X,Y)$, we deduce
		\begin{align*}
		& \| U_{\varepsilon}(t,\tau) - U_0(t,\tau) \|_{\mathcal{L}(X)} \\
		%
		%
		%Primeira Desigualdade
		%
		%
		& \quad \leq  \| T_{-A_{\varepsilon} (\tau) } (t-\tau) - T_{-A_0(\tau)} (t-\tau)  \|_{\mathcal{L}(X,Y)}
		+
		\int_{\tau}^{t} 
		\| 
		T_{-A_{\varepsilon}(s)}(t-s) 	\|_{\mathcal{L}(X,Y)}
		\|
		[ \Phi_{\varepsilon} (s,\tau)-\Phi_0(s,\tau)]
		\|_{\mathcal{L}(X)}
		ds \\
		& \qquad +
		\int_{\tau}^{t}
		\| 
		[T_{-A_{\varepsilon}(s)} (t-s) - T_{A_{0} (s)} (t-s)] 
		\|_{\mathcal{L}(X,Y)}
		\|
		\Phi_0(s,\tau)
		\|_{\mathcal{L}(X)}
		ds\\
		%
		%
		%Segunda Desigualdade
		%
		%
		& \quad \leq C(t-\tau)^{-1+\beta(1-\theta)} [\eta(\varepsilon)]^{\theta}
		+
		\int_{\tau}^{t} 
		C
		(t-s)^{\beta-1} (s-\tau)^{-1+\delta(1-\theta)} e^{K(s-\tau)} \ell (\theta, \varepsilon)
		ds \\
		& \quad +
		\int_{\tau}^{t}
		C(t-s)^{-1+\beta(1-\theta)} [\eta(\varepsilon)]^{\theta}
		(s-\tau)^{\delta-1}e^{K(s-\tau)}
		ds\\
		%
		%
		%Terceira Desigualdade
		%
		%
		& \quad \leq C(t-\tau)^{-1+\beta(1-\theta)} [\eta(\varepsilon)]^{\theta}
		+
		C e^{K(t-\tau)} \ell (\theta, \varepsilon)
		(t-\tau)^{\beta+\delta(1-\theta)-1} \mathcal{B}(\beta,\delta(1-\theta)) \\
		& \qquad+
		C [\eta(\varepsilon)]^{\theta}
		(t-\tau)^{\delta + \beta(1-\theta)-1} \mathcal{B} (\beta(1-\theta), \delta) e^{K(t-\tau)} \\
		%
		%
		%Quarta Desigualdade
		%
		%
		& \quad \leq C(\theta, \delta) (t-\tau)^{-1+\beta(1-\theta)} e^{K(t-\tau)} \ell(\theta, \varepsilon).
		\end{align*}
	{	\hspace{16.5cm}	\qedsymbol}
%	\end{proof}
%\end{theorem}

\subsection{Convergence and estimates of the solution of the semilinear problem}

As mentioned at Section \ref{S:Setting_and_main_results}, Problem \eqref{Pe} is locally well-posed, that is, there exists 
$$u_{\varepsilon}(t,\tau,u^{\tau}) = U_{\varepsilon}(t,\tau) u^{\tau} + \int_{\tau}^{t} U_{\varepsilon}(t,s) F_{\varepsilon}(s, u_{\varepsilon}(s)) ds,$$ 
that solves the problem for $t \in [\tau, \tau+T(\varepsilon,\tau, u^{\tau}))$. We actually have global well-posedness ($T(\varepsilon,\tau,u^{\tau}) = \infty$) as a consequence of the following result.

\begin{lemma}\label{L_est_sol}
	Let $\beta \in (0,1]$ be the constant in the resolvent estimate \eqref{res_est_X,Y}. There exist constants $C,K>0$ such that, for any $\tau \in \mathbb{R}$ and $T>0$ for which $u_{\varepsilon}(\cdot, \tau, u^{\tau})$ is defined in  $ (\tau, \tau+T]$, we have
	\begin{equation*}\label{est_sol}
	\| u_{\varepsilon}(t, \tau, u^{\tau}) \|_{Y }
	\leq 
	C(t-\tau)^{\beta-1}e^{K(t-\tau)}[1+\|u^{\tau}\|_Y], \qquad \mbox{for all } t\in (\tau, \tau+T]. 
	\end{equation*} 
	
	In particular, the $\|u_{\varepsilon}(t)\|_Y$ does not blow-up in any finite time interval and it is globally defined.
	
	\begin{proof}
	This result follows from the expression for the solution $u_{\varepsilon}$ and the estimates obtained previously.
	\begin{align*}
	\|u_{\varepsilon}(t,\tau,u^{\tau})\|_Y 
	& \leq \|U_{\varepsilon}(t,\tau)\|_{\mathcal{L}(X,Y)}\|u^{\tau}\|_X + \int_{\tau}^{t}\|U_{\varepsilon}(t,s)\|_{\mathcal{L}(X,Y)}\|F_{\varepsilon}(s,u_{\varepsilon}(s,\tau,u^{\tau}))) \|_X ds\\
	& \leq C(t-\tau)^{\beta-1}e^{K(t-\tau)} \|u^{\tau}\|_Y+\int_{\tau}^{t} C(t-s)^{\beta-1}e^{K(t-s)} M ds \\
	&\leq C(t-\tau)^{\beta-1} e^{K(t-\tau)}\|u^{\tau}\|_Y + \frac{CM}{\beta} (t-\tau)^{\beta} e^{K(t-\tau)}\\
	& \leq C(t-\tau)^{\beta-1} e^{K(t-\tau)} [1+\|u^{\tau}\|_Y].
 	\end{align*}
 	
 	Therefore, the solution is bounded in any bounded interval $[\tau+m, \tau+T]$, for $0<m<T$, being globally defined.
	\end{proof}
\end{lemma}

We now prove Theorem \ref{T_Conv_solution} that provides a rate at which the solutions converge.

\

\noindent \textbf{Proof of Theorem \ref{T_Conv_solution}:} In the sequel we will denote the solution $u_{\varepsilon}(t,\tau,u^{\tau})$  by  $u_{\varepsilon}(t)$. Let $M>0, L>0$ be the boundedness and Lipschitz constant for $F$, respectively, and 
$\rho(\theta,\varepsilon) = \max\{ [\eta(\varepsilon)]^{\theta}, [\xi(\varepsilon)]^{\theta}, \gamma(\varepsilon)  \}.$
 Using the expression for the solution and rates of convergence established earlier, we obtain
\begin{align*}
\| u_{\varepsilon}(t) -u_{0}(t) \|_Y 
%
%Primeira igualdade
%
& 
\leq \|U_{\varepsilon}(t,\tau) - U_{0}(t,\tau)\|_{\mathcal{L}(X,Y)} \|u^{\tau}\|_X + \int_{\tau}^{t} \|U_{\varepsilon}(t,s) - U_0(t,s)\|_{\mathcal{L}(X,Y)}  \|F_0 (s,u_0(s)) \|_X ds \\
& \quad + \int_{\tau}^{t} \|U_{\varepsilon}(t,s)\|_{\mathcal{L}(X,Y)} \left[ \| F_{\varepsilon} (s,u_{\varepsilon}(s)) - F_{\varepsilon}(s,u_0(s))\|_X + \|F_{\varepsilon}(s,u_0(s)) - F_0 (s,u_0(s)) \|_X \right]ds \\
%
%
%Segunda Desigualdade
%
& 
\leq C(t-\tau)^{-1+\beta (1-\theta)} e^{K(t-\tau)} \ell(\theta,\varepsilon) \|u^{\tau}\|_Y + \int_{\tau}^{t} C(t-s)^{-1+\beta(1-\theta)}e^{K(t-s)} \ell(\theta,\varepsilon) M ds\\
& \quad +  \int_{\tau}^{t} C (t-s)^{\beta-1}e^{K(t-s)}\left[ L\|  u_{\varepsilon}(s)- u_0(s)\|_Y + \gamma(\varepsilon) \right]ds \\
%
%Terceira desigualdade
%
& \leq C(t-\tau)^{-1+\beta (1-\theta)} e^{K(t-\tau)} \ell(\theta,\varepsilon) \|u^{\tau}\|_Y +  \frac{CM}{\beta(1-\theta)}(t-\tau)^{\beta(1-\theta)}e^{K(t-s)} \ell(\theta,\varepsilon) \\
& \quad + \frac{C}{\beta}e^{K(t-\tau)} (t-\tau)^{\beta} \gamma(\varepsilon) + CL\int_{\tau}^{t}  (t-s)^{\beta-1}e^{K(t-s)} \|  u_{\varepsilon}(s)- u_0(s)\|_Y ds \\
%
%Quarta Desigualdade
%
& \leq C(t-\tau)^{-1+\beta (1-\theta)} e^{K(t-\tau)} 
\left[ \|u^{\tau}\|_Y +1 \right] \rho(\theta, \varepsilon) 
 + CL\int_{\tau}^{t}  (t-s)^{\beta-1}e^{K(t-s)} \|  u_{\varepsilon}(s)- u_0(s)\|_Y ds, \\
\end{align*}
where we incorporated the terms $(t-\tau)$ with a positive exponent to the exponential growth given by $e^{K(t-\tau)}$, making adjustments in the constant $C$, if necessary. Multiplying both sides by $e^{-K(t-\tau)}$ and considering $\Psi (t) = e^{-K(t-\tau)} \| u_{\varepsilon}(t) - u_0 (t) \|_Y$, we obtain
	\[
	\Psi(t) \leq C(t-\tau)^{-1+\beta(1-\theta)} [\|u^{\tau}\|_Y +1 ] \rho (\theta, \varepsilon) +CL \int_{\tau}^{t} (t-s)^{\beta-1} \Psi(s) ds.
	\]

	We now apply Gronwall's inequality \cite[p.190]{henry} to conclude that 
	\[
	\Psi(t) \leq \frac{C}{\beta (1-\beta)}(t-\tau)^{-1+\beta(1-\theta)} [\|u^{\tau}\|_Y +1 ] \rho (\theta, \varepsilon)e^{\tilde{K}(t-\tau)},
	\]
	where $\tilde{K}>(2CL\Gamma(\beta))^{\frac{1}{\beta}}$. Therefore,
	\[
	\| u_{\varepsilon}(t) - u_0(t)\|_Y  \leq \frac{C}{\beta (1-\beta)}(t-\tau)^{-1+\beta(1-\theta)} [\|u^{\tau}\|_Y +1 ] \rho (\theta, \varepsilon)e^{(\tilde{K} +K)(t-\tau)}
	\]
	and 
	\[
		\| u_{\varepsilon}(t) - u_0(t)\|_Y \stackrel{\varepsilon \to 0}{\longrightarrow} 0,
	\]
	uniformly for $t$ in compact subsets of $(\tau, \infty)$, any $\tau \in \mathbb{R}$ and $u^{\tau}$ in bounded sets of $Y$.

	\hspace{16.65cm} \qedsymbol

%%%%%%%%%%%%%%%%%%%%%%%%%%%%%%%%%%%%%%%%%%%%%%%%%%%%
%SECTION 3 - REACTION-DIFFUSION EQUATION
%%%%%%%%%%%%%%%%%%%%%%%%%%%%%%%%%%%%%%%%%%%%%%%%%%%%

\section{Application to  reaction-diffusion equations with varying diffusion coefficients}\label{S:Case_I_application_I}

	As a first application of the abstract theory developed in the previous sections, we consider a family in $\varepsilon \in [0,1]$  of  singularly nonautonomous reaction-diffusion equation in a bounded smooth domain $\Omega \subset \mathbb{R}^{3}$
	\begin{equation}\label{eq_1*}
	\begin{split}
	& (u_{\varepsilon})_t - div (a_{\varepsilon} (t,x) \nabla u_{\varepsilon}) +u_{\varepsilon} = f_{\varepsilon} (t,u_{\varepsilon}), \hspace{1.1cm} x\in \Omega,  t> \tau,\\
	& \partial_n u_{\varepsilon} =0,  \hspace{5.9cm} x\in \partial \Omega,\\
	& u_{\varepsilon} (\tau,x) = u^{\tau} (x). 
	\end{split}
	\end{equation}
	
	An autonomous version (where $a_{\varepsilon}$ and $f_{\varepsilon}$ do not depend on $t$) was completely studied in \cite{ArrietaBezerraCarvalho_2013} and the authors obtained rate of convergence of solutions and attractors in terms of $\varepsilon$. The nonautonomous counterpart \eqref{eq_1*} was introduced in \cite{BCNS_2}, where the authors studied global well-posedness and existence of pullback attractor, but for a single equations rather than a family of equations parametrized in $\varepsilon \in [0,1]$.

	We shall apply the abstract theory developed in Section \ref{S:Setting_and_main_results} in order to obtain a rate at which solutions of \eqref{eq_1*} converge as ${\varepsilon \to 0^{+}}$. We assume the following conditions for the problem:

\begin{enumerate}[label={(A.\arabic*)},ref=(A.\arabic*)]

	\item \label{A_1} The functions	$a_{\varepsilon} :\R \times \overline{ \Omega  }  \rightarrow \R^{+}$ are continuously differentiable with respect to the second variable, and $a_{\varepsilon}(\cdot,\cdot)$ has its image in a closed interval $[m,M] \subset (0,\infty)$. We also assume that the gradient function (in $x$) of $a_{\varepsilon}(t,x)$ is bounded, that is, $ \nabla_x a_{\varepsilon}(t,x) \in [L^{\infty} ({ \Omega  })]^{3}$.

	\item \label{A_2} Both functions $a_{\varepsilon}(\cdot, \cdot)$ and $\nabla_x a_{\varepsilon} (\cdot,\cdot)$ are uniformly $\delta-$\emph{H\"{o}lder continuous in the first variable} that is, there exists $\delta\in (0,1]$ and a constant $C>0$ such that
	\begin{equation*}\label{HC_a_b}
	|a_{\varepsilon}(t,x) - a_{\varepsilon}(s,x)| \leq C|t-s|^{\delta}, \quad |\nabla_xa_{\varepsilon}(t,x) - \nabla_xa_{\varepsilon}(s,x)| \leq C|t-s|^{\delta}, 
	\end{equation*}
	for all $\varepsilon \in[0,1]$, $t,s\in \mathbb{R}$ and $x\in \Omega$.

	\item \label{A_3}  For each $\varepsilon \in [0,1]$, $f_\varepsilon \in \mathcal{C}^{1}(\mathbb{R}\times \mathbb{R}, \mathbb{R})$ and satisfies a \emph{polynomial growth condition of order $\rho$}, that is, there exists $C$ and  $1\leq \rho <3$ such that 
	
	\vspace{-0.4cm}
	\begin{align*}
	|f_{\varepsilon}(t,\xi)-f_{\varepsilon}(t,\psi)| &\leq C|\xi-\psi|(1+|\xi|^{\rho -1}+|\psi|^{\rho -1}) , \\
	|f_{\varepsilon}(t,\xi)| &\leq C  (1+|\xi|^{\rho}).
	\end{align*}
	
	%\noindent In particular, it follows from $f_{\varepsilon} \in C^{1}(\mathbb{R}\times\mathbb{R}, \mathbb{R})$  that  $t\mapsto f_{\varepsilon}(t, \cdot )$ is  \emph{locally H\"{o}lder continuous}.

		\item \label{A_4} We define the quantities
		%		
		%ESTIMATIVA PARA a_e	
		%
		$$
		\|a_{\varepsilon} - a_{0}\|_{\infty} := 
		\sup_{t \in \mathbb{R}}
		\| a_{\varepsilon}(t,\cdot) - a_0(t,\cdot)\|_{L^{\infty} (\Omega)}, $$
		%
		%Estimativa para \nabla_x a_e
		%
		$$
		\| \nabla_x a_{\varepsilon} - \nabla_x a_0\|_{\infty} := 
		\sup_{t \in \mathbb{R}}
		\| \nabla_x a_{\varepsilon}(t,\cdot) - \nabla_x a_0(t,\cdot)\|_{[L^{\infty} (\Omega)]^{3}}, $$
		%
		%Estimativa para a f_e
		%
		$$\| f_{\varepsilon} - f_0\|_{\infty}
		:= 
		\sup_{t\in \mathbb{R} } \| f_{\varepsilon}(t,\cdot) - f_0(t,\cdot)\|_{L^{\infty}(\Omega) },$$ 
	
		\noindent and we assume that each one of them varies continuously on $\varepsilon \in [0,1]$. In particular,they approach to zero as $\varepsilon \to 0^{+}$. 
	
\end{enumerate}

	The upper bound requested for $\rho $ in \ref{A_3} will become clear after we specify the  phase space in which we pose the problem.  Under conditions above, we write Problem \eqref{eq_1*} in its abstract form as follows: the linear part of the equation (which is time-dependent) is given by the operator
	$A_{\varepsilon}(t):D(A_{\varepsilon}(t)) \subset L^{2}(\Omega) \rightarrow L^{2}(\Omega)$ where 
	\begin{align*}
	&D=D(A_{\varepsilon}  (t) )=  \left\{ u \in H^{2}(\Omega) :
	\partial_n u = 0 \mbox{ in }\partial \Omega \right\} =: H^{2}_{\mathcal{N}}
	,
	%\label{D_A_e} 
	\\
	&\hspace{0.2cm} A_{\varepsilon}(t) u  =  -div( a_{\varepsilon}(t,x)\nabla u )+ u, \mbox{ for }u \in D. %\label{A_e}
	\end{align*}

This family of linear operators has well-known properties that we gather in the sequel. They follow from classical spectral theory (see \cite{Arrieta_et_al_1999,Kato,pazy}) and from the properties required upon $a_{\varepsilon}$ in \ref{A_1} and \ref{A_2}. To simplify notation, we shall omit the domain $\Omega$ in the space norms, that is, $\|\cdot\|_{L^{2}} = \|\cdot\|_{L^{2}(\Omega)}$.  

\begin{proposition}\label{Prop_A}
	This family $\{ A_{\varepsilon}(t), \ t\in \mathbb{R}\}_{\varepsilon \in [0,1]}$ has the following properties:
	\begin{enumerate}
		\item $D(A_{\varepsilon}(t))$ does not depend on $t$ or $\varepsilon$. Moreover, for any fixed $\varepsilon \in [0,1]$ and $t \in \mathbb{R}$,  the graph norm $\| A_{\varepsilon}(t) \cdot \|_{L^{2}}$ is equivalent to $H^{2}(\Omega)-$norm when restricted to $D$, that is, for any $u \in D$,
		\begin{equation*}
		C_1 \|u \|_{H^{2}} \leq \| A_{\varepsilon} (t) u \|_{L^{2}} \leq  C_2 \| u \|_{H^{2}},
		\end{equation*}
		and constants $C_1,C_2$ are uniform for $\varepsilon \in [0,1]$ and $t \in \mathbb{R}$.
		\item $A_{\varepsilon} (t)$ is self-adjoint 	and has compact resolvent.
		\item \label{4} Its spectrum consists entirely of isolated eigenvalues, all of them positive and real, with the first being $1$:
		$$\sigma(A_{\varepsilon}(t)) =\{  \lambda_{\varepsilon, i}(t); \ i \in \mathbb{N}^{*} \mbox{ and } 1 = \lambda_{\varepsilon,1} (t) \leq \lambda_{\varepsilon,2}(t) \leq ... \leq \lambda_{\varepsilon,n}(t) \leq ...   \} .$$
		
		\item \label{5} 
		For any $\frac{\pi}{2}<\varphi < \pi$,  $\Sigma_{\varphi} = \left\{  \lambda \in \mathbb{C}; |\arg \lambda| \leq \varphi \right\} \subset \rho(-A_{\varepsilon}(t))$ and 
		\begin{align*}
		 \| (\lambda I + A_{\varepsilon}(t) )^{-1}\|_{\mathcal{L}(L^{2})} & \leq \frac{C}{|\lambda|^{} +1} ,\quad \forall \lambda\in \Sigma_{\varphi} \cup\{0\}, %\label{lim_res_L2} 
		 \\
		\| (\lambda I + A_{\varepsilon}(t) )^{-1}\|_{\mathcal{L}(H^{1})}  & \leq \frac{C}{|\lambda|^{} +1},\quad \forall \lambda\in \Sigma_{\varphi} \cup\{0\}, %\label{lim_res_H1} 
		\\
		\| (\lambda I + A_{\varepsilon}(t) )^{-1}\|_{\mathcal{L}(L^{2},H^{1})}  & \leq \frac{C}{|\lambda|^{\frac12} +1}, \quad \forall \lambda\in \Sigma_{\varphi} \cup\{0\}, %\label{lim_res_L2_H_1}
		\end{align*}
		where $C$ does not depend on $\varepsilon$ or $t$ (only on $\varphi$).
	\end{enumerate}

\end{proposition}

We restate Problem \eqref{eq_1*} as an  abstract  semilinear evolution problem:
\vspace{0.2cm}
\begin{equation}\label{eq_1_abs}
\begin{split}
&(u_{\varepsilon})_t + A_{\varepsilon}(t) u_{\varepsilon} = F_{\varepsilon}(t,u_{\varepsilon}), \quad  t> \tau,\\
&u_{\varepsilon}(\tau) = u^{\tau} \in H^{1}(\Omega),  
\end{split}
\end{equation}

\vspace{0.1cm}\noindent where $F_{\varepsilon}$ is a nonlinearity given by 
$$F_{\varepsilon}(t,u_{\varepsilon}) (x) = f_{\varepsilon}(t,u_{\varepsilon}(t,x)).$$ 

Since 
$$H^{1}(\Omega) \hookrightarrow L^{r}(\Omega), \mbox{ for all } 2\leq r < 6 
,$$ 

\noindent then the growth condition \ref{A_3} required for $f_{\varepsilon}$ implies that $F:\mathbb{R} \times H^{1}(\Omega) \rightarrow L^{2}(\Omega)$, as long as $1\leq \rho < 3$. With the notation of Section \ref{S:Setting_and_main_results}, $L^2(\Omega)$ will play the role of Banach space $X$ and $H^{1}(\Omega)$ the Banach space $Y$. Moreover, one can easily check that from \ref{A_3} we derive, for any $\varepsilon\in [0,1]$ and $t\in \mathbb{R}$,
\begin{align*}
\|F_{\varepsilon}(t,u) - F_{\varepsilon}(t,v) \|_{L^{2}} 
& \leq C \| u-v \|_{H^{1}} \left[ 1+ \|u\|_{H^{1}}^{\rho-1} + \| v \|^{\rho-1}_{H^{1}} \right], \\
\|F_{\varepsilon}(t,u) \|_{L^{2}} & \leq C\left[  1+ \|u\|^{\rho}_{H^{1}}  \right].
\end{align*}

In order to apply the theory developed in Section \ref{S:Setting_and_main_results}, we first need to verify that \ref{P_1} to\ref{P_5} hold for \eqref{eq_1_abs}. From Proposition \ref{Prop_A}, properties 
\ref{P_1} and \ref{P_2} already follow. Property \ref{P_3} is proved in next lemma.

\begin{lemma}\label{HC_A_0}
	Assume \ref{A_1} and \ref{A_2} hold and let $\delta \in (0,1]$ be the uniform H\"{o}lder exponent for $t\mapsto a_{\varepsilon}(t,\cdot)$ and $t\mapsto \nabla_x a_{\varepsilon}(t,\cdot)$. Then, there exists a constant $C>0$, independent of $\varepsilon \in [0,1]$ or $\tau \in \mathbb{R}$, such that,  for any $\varepsilon_1, \varepsilon_2 \in [0,1]$, the function 
	$$\mathbb{R} \ni t \mapsto A_{\varepsilon_1}(t)A_{\varepsilon_2}(\tau)^{-1} \in \mathcal{L} (L^{2}(\Omega))$$
	is H\"{o}lder continuous with exponent $\delta$, that is, 
	\[
	\| [A_{\varepsilon_1}(t) - A_{\varepsilon_1}(s)]  A_{\varepsilon_2}(\tau)^{-1} \|_{\mathcal{L}(L^{2})} \leq C |t-s|^{\delta}, \quad \mbox{for all }  \tau, s, t\in \mathbb{R}.
	\]

	\begin{proof}
		
		For any $u\in D$, we have
		$
		A_{\varepsilon_1}(t)u - A_{\varepsilon_1}(s)u = 
		-div \left( [a_{\varepsilon_1}(t,x)-a_{\varepsilon_1}(s,x)] \nabla u \right)
		$
		and
		
		\begin{equation*}
		\begin{split}
		\| A_{\varepsilon_1}(t) u - A_{\varepsilon_1}(s) u \|_{L^{2}}^{2} & = \int_{\Omega}
		\left| 
		div([a_{\varepsilon_1}(t,x)  - a_{\varepsilon_1}(s,x)]\nabla u(x))   
		\right|^{2}
		dx \\
		%
		% Segunda igualdade
		%
		&  = 
		\int_{\Omega}
		\left| 
		\nabla_x ([a_{\varepsilon_1}(t,x)  - a_{\varepsilon_1}(s,x)]) \nabla u(x) +    [a_{\varepsilon_1}(t,x)  - a_{\varepsilon_1}(s,x)] \Delta u (x)
		\right|^{2}
		dx \\
		%
		%Terceira desigualdade
		%
		&  \leq C |t-s|^{2\delta }
		\int_{\Omega}
		\left\{
		\frac{\left| \nabla_x a_{\varepsilon_1}(t,x)  - \nabla_x a_{\varepsilon_1}(s,x)  \right| }{ |t -s|^{\delta}}
		\right\}^{2}
		|\nabla u(x) |^{2} 
		dx \\
		&  \qquad +C |t-s|^{2\delta }
		\int_{\Omega}
		\left\{
		\frac{\left| a_{\varepsilon_1}(t,x)  -  a_{\varepsilon_1}(s,x)  \right| }{ |t -s|^{\delta}}
		\right\}^{2}
		|\Delta u(x) |^{2} 
		dx \\
		%
		%Quarta desigualdade
		%
		&  \leq C |t-s|^{2\delta } 
		\left\{  
		\| \nabla u \|_{L^{2}}^{2}
		+\| \Delta u \|_{L^{2}}^{2}
		\right\}
		%
		%
		%
		%& \qquad 
		\leq C|t-s|^{2\delta } \| u \|_{H^{2}}^{2}. 
		\end{split}
		\end{equation*}

		Taking the square root on both sides and replacing $u$ by $A_{\varepsilon_2}(\tau)^{-1}w$, then we have for any $w\in L^{2}(\Omega)$

		\begin{equation*}
		\|[{A_{\varepsilon_1}}(t) - {A_{\varepsilon_1}}(s)]A_{\varepsilon_2}(\tau)^{-1}w\|_{L^{2}} 
		\leq
		C|t-s|^{ \delta} \|A_{\varepsilon_2} (\tau)^{-1} w\|_{H^{2}} 
		\leq 
		C|t-s|^{ \delta} \|w\|_{L^{2}}.
		%\quad \forall w\in L^{2}(\Omega).
		\end{equation*}
	\end{proof}

\end{lemma}

	It remains to check the properties responsible to make the connections among the problems as $\varepsilon$ varies in $[0,1]$. Those are conditions \ref{P_4} and \ref{P_5} from Section \ref{S:Setting_and_main_results}. We start verifying \ref{P_4} and we begin by proving an auxiliary result.

	\begin{lemma}\label{L_P_4_pt1}
	Let $\varepsilon_1, \varepsilon_2 \in [0,1]$ and $t,\tau \in \mathbb{R}$. Then
	\begin{equation*}
	\| [A_{\varepsilon_1}(t)  - A_0(t)  ] A_{\varepsilon_2}(\tau)^{-1}\|_{\mathcal{L}(L^{2})} 
	\leq
	C ( \| a_{\varepsilon_1}-a_0 \|_{\infty} + \| \nabla_x a_{\varepsilon_1} - \nabla_x a_0 \|_{\infty}).
	\end{equation*}

	\begin{proof}
	Take $u = A_{\varepsilon_2} (\tau)^{-1} w \in D$, where  $w$ is any element in $ L^{2}(\Omega)$. We have
	\begin{align*}
	\| [A_{\varepsilon_1}(t) - A_0(t)] u \|_{L^{2}}^{2} 
	&
	\leq \int_{\Omega} | div [(a_{\varepsilon_1} (t,x)-a_0(t,x) ) \nabla u ]|^{2} dx \\
	&\leq  C \int_{\Omega}  [|\nabla_x (a_{\varepsilon_1}(t,x) - a_0(t,x) ) |^{2} | \nabla u|^{2} + | a_{\varepsilon_1}(t,x) - a_0(t,x)|^{2} |\Delta u|^{2} ]
	dx\\
	& \leq C ( \| a_{\varepsilon_1} - a_0 \|^{2}_{\infty} + \| \nabla_x a_{\varepsilon_1} -\nabla_x a_{0} \|_{\infty}^{2} ) \|u\|_{H^{2}}^{2}.
	\end{align*}
	
	Therefore,
	\begin{align*}
	\| [A_{\varepsilon_1}(t) - A_0(t)] A_{\varepsilon_2}(\tau)^{-1} w \|_{L^{2}} 
	& \leq C ( \| a_{\varepsilon_1} - a_0 \|_{\infty} + \| \nabla_x a_{\varepsilon_1} -\nabla_x a_{0} \|_{\infty} ) \|A_{\varepsilon_2}(\tau)^{-1} w\|_{L^{2}}\\
	& \leq C ( \| a_{\varepsilon_1} - a_0 \|_{\infty} + \| \nabla_x a_{\varepsilon_1} -\nabla_x a_{0} \|_{\infty} ) \| w\|_{L^{2}}.
	\end{align*}
	\end{proof}
	\end{lemma}
	
	With the previous lemma, we are now able to prove that \ref{P_4} holds for this problem.

	\begin{lemma}\label{P_P_4}
	Let $t,\tau \in \mathbb{R}$, $\varepsilon \in [0,1]$. There exists $C>0$ independent of $t,\tau,\varepsilon$ such that
	\begin{equation*}
	\| A_{\varepsilon}(t) A_{\varepsilon}(\tau)^{-1}-A_0(t)A_0(\tau)^{-1} \|_{\mathcal{L}{(L^{2})} }\leq C( \| a_{\varepsilon}-a_0 \|_{\infty} + \| \nabla_x a_{\varepsilon} - \nabla_x a_0 \|_{\infty}).
	\end{equation*}
	
	\begin{proof}
	Let $w\in L^{2}(\Omega)$ and consider $u = A_{\varepsilon}(\tau)^{-1} w$, $v= A_0(\tau)^{-1}w$. It follows from the boundedness of $a_{\varepsilon}$ and its convergence to $a_0$ that
	\begin{align*}
	\|A_{\varepsilon}(t)u - A_0(t)v \|_{L^{2}}^{2} 
	& \leq C \int_{\Omega} |
	 div [ a_{\varepsilon}(t,x) \nabla u  ] - div [ a_0(t,x) \nabla v ]|^{2} dx + C \int_{\Omega} |u-v|^{2} dx \\
	 & \leq C \int_{\Omega} | div[ (a_{\varepsilon}(t,x) -a_0(t,x)) \nabla u] + div [a_0(t,x) (\nabla u - \nabla v) ]   |^{2}dx +
	  C \int_{\Omega} |u-v|^{2} dx \\ 
	  & \leq C \int_{\Omega} \scalebox{1.5}{$\{$}  | \nabla_x a_{\varepsilon}(t,x) - \nabla_x a_0(t,x) |^{2}
	  |\nabla u|^{2} + |a_{\varepsilon}(t,x)-a_0(t,x)|^{2} |\Delta u|^{2} \\
	  & \hspace{1.75cm} +|\nabla_x a_0(t,x) |^{2} | \nabla (u-v)|^{2} + |a_0(t,x)|^{2} |\Delta (u-v)|^{2} \scalebox{1.5}{$\}$}dx
	  +
	  C \int_{\Omega} |u-v|^{2} dx \\ 
	  & \leq C \| \nabla_x a_{\varepsilon} - \nabla_x a_0 \|_{\infty}^{2} \left( \int_{\Omega} |\nabla u|^{2}dx  \right) + C\| a_{\varepsilon} -a_0 \|_{\infty}^{2} \left( \int_{\Omega} |\Delta u|^{2}  \right) + C \| u-v\|_{H^{2}}^{2} \\
	  &\leq C ( \|a_{\varepsilon} - a_{0}\|_{\infty}^{2}+\|\nabla_xa_{\varepsilon}-\nabla_x a_0 \|_{\infty}^{2})\|u\|^{2}_{H^{2}}+C \|u-v\|_{H^{2}}^{2}.
	\end{align*}

	From the choice of $u,v$ and from Lemma \ref{L_P_4_pt1} we obtain
	\begin{align*}
	& \| A_{\varepsilon}(t) A_{\varepsilon}(\tau)^{-1}w - A_0(t) A_0(\tau)^{-1} w \|_{L^{2}}^{2}  \\
	& \quad \leq C ( \|a_{\varepsilon} - a_{0}\|_{\infty}^{2}+\|\nabla_xa_{\varepsilon}-\nabla_x a_0 \|_{\infty}^{2})\|A_{\varepsilon}(\tau)^{-1} w\|^{2}_{H^{2}}
	+C \|A_{\varepsilon}(\tau)^{-1}w - A_0(\tau)^{-1} w\|_{H^{2}}^{2}\\
	& \quad \leq C ( \|a_{\varepsilon} - a_{0}\|_{\infty}^{2}+\|\nabla_xa_{\varepsilon}-\nabla_x a_0 \|_{\infty}^{2})\| w\|^{2}_{H^{2}}
	+C \|[A_0(\tau) - A_{\varepsilon}(t)]A_0(\tau)^{-1} w\|_{H^{2}}^{2}\\
	& \quad \leq 2C ( \|a_{\varepsilon} - a_{0}\|_{\infty}^{2}+\|\nabla_xa_{\varepsilon}-\nabla_x a_0 \|_{\infty}^{2})\| w\|^{2}_{H^{2}}.
	\end{align*}
	\end{proof}
	
	\end{lemma}
	
	Therefore, Lemma \ref{P_P_4} states that \ref{P_4} holds for 
	$$ \xi(\varepsilon)=C( \| a_{\varepsilon_1}-a_0 \|_{\infty} + \| \nabla_x a_{\varepsilon} - \nabla_x a_0 \|_{\infty}).$$

	 Inspired in \cite{ArrietaBezerraCarvalho_2013}, we will use a variational formulation in order to obtain resolvent convergence
	\begin{equation*}%\label{conv_inverso_ex_1}
	\|  A(t)^{-1}_{\varepsilon} - A(t)^{-1}_0 \|_{ \mathcal{L}(L^{2},H^{1})} \stackrel{ \varepsilon \to 0^{+}  }{\longrightarrow} 0,
	\end{equation*}
	that is, in order to prove that \ref{P_5} holds.

\begin{lemma}
	Given $ g \in X=L^{2}(\Omega)$, a fixed $t\in \mathbb{R}$ and $\varepsilon \in [0,1]$, there exists a unique $u_{\varepsilon} \in  H^{2}_{\mathcal{N} }$ solution of 
	\begin{equation}\label{eliptic_problem_ex_1}
	\begin{cases}
	- div (a_{\varepsilon} (t,x) \nabla u_{\varepsilon} ) + u_{\varepsilon} = g, \quad x\in \Omega, \\
	\partial_n u_{\varepsilon} = 0, \quad \partial \Omega.
	\end{cases}
	\end{equation}
	
	Moreover,
	\begin{enumerate}
		\item there exists $C>0$, independent of $\varepsilon \in [0,1]$, $g\in L^{2}(\Omega)$ and $t\in \mathbb{R}$, such that
		$$\| u_{\varepsilon} \|_{H^{1}} \leq C \| g\|_{L^{2}}.
		$$
		
		\item There is also a constant $C>0$, independent of $\varepsilon\in [0,1]$, $g\in L^{2}(\Omega)$  and $t\in \mathbb{R}$, such that 
		\begin{equation*}\label{conv_1}
		\| u_{\varepsilon} - u_{0}\|_{H^{1}} \leq C \|a_{\varepsilon} - a_0 \|_{\infty} \|g\|_{L^{2}}.
		\end{equation*}
		
	\end{enumerate}

	\begin{proof}
		Existence of $u_{\varepsilon}$ that solves \eqref{eliptic_problem_ex_1} follows from the fact that $0\in \rho (A_{\varepsilon}(t))$, for all $t\in \mathbb{R}$ and $\varepsilon \in [0,1]$. To prove the first statement, we consider the weak formulation of \eqref{eliptic_problem_ex_1}:
		\begin{equation*}
		\int_{\Omega} a_{\varepsilon} (t,x) \nabla u_{\varepsilon} \nabla \varphi + \int_{\Omega} u_{\varepsilon} \varphi = \int_{\Omega} g \varphi, \quad \mbox{ for } \varphi \in H^{1}(\Omega).
		\end{equation*}
		
		By taking $\varphi = u_{\varepsilon}$, using Young's inequality and the fact that $a_{\varepsilon}(\cdot, \cdot) \subset [m, M]$, we obtain, for any $\nu>0$,
		\begin{align*}
		& \int_{\Omega} a_{\varepsilon} (t,x) [\nabla u_{\varepsilon}]^{2} 
		+ \int_{\Omega} [u_{\varepsilon}]^{2} 
		= \int_{\Omega} g u_{\varepsilon}, \\
		& m \int_{\Omega} [\nabla u_{\varepsilon}]^{2} 
		+ \int_{\Omega} [u_{\varepsilon}]^{2} 
		\leq  \int_{\Omega} |g| |u_{\varepsilon}| 
		\leq \int_{\Omega} \left[ \left\{    \frac{1}{\nu^{2}} \frac{|g|^{2}}{2}   \right\} + \left\{  \frac{\nu^{2} |u_{\varepsilon}|^{2} }{2}  \right\}   \right],\\
		& m \int_{\Omega} [\nabla u_{\varepsilon}]^{2} 
		+ \left(1 -   \frac{\nu^{2} }{2} \right) \int_{\Omega} [u_{\varepsilon}]^{2} 
		\leq   \frac{1}{2\nu^{2}}  \int_{\Omega}  |g|^{2}. \\
		\end{align*}
		
		Choosing $\nu$ small such that $1 - \frac{\nu^{2} }{2} >0$  we obtain 
		\begin{equation}\label{uniform_bound}
		 \|u_{\varepsilon}\|^{2}_{H^{1}} \leq C \|g\|^{2}_{L^{2}}, 
		\end{equation}
		where $C = \frac{1 }{ 2\nu^{2}  \min\{ m,1 - \frac{\nu^{2}}{2} \} } $, which does not depend on $\varepsilon$ or $t$. 
		
		For the second statement we proceed similarly. Rather than taking $u_{\varepsilon}$ as a test function, we choose $u_{\varepsilon} - u_{0}$, obtaining 
		
		\begin{align*}
		&\int_{\Omega} a_{\varepsilon} (t,x) \nabla u_{\varepsilon} \left(  \nabla u_{\varepsilon} - \nabla u_{0} \right) + \int_{\Omega} u_{\varepsilon} \left(u_{\varepsilon} - u_{0} \right) = \int_{\Omega} g \left(  u_{\varepsilon} - u_{0} \right),\\
		&	\int_{\Omega} a_{0} (t,x) \nabla u_{0} \left(  \nabla u_{\varepsilon} - \nabla u_{0} \right) + \int_{\Omega} u_{0} \left(u_{\varepsilon} - u_{0} \right) = \int_{\Omega} g \left(  u_{\varepsilon} - u_{0} \right).
		\end{align*}

		Equality on the right side implies 
		\begin{align*}
		\int_{\Omega} a_{\varepsilon} (t,x) \nabla u_{\varepsilon} \left(  \nabla u_{\varepsilon} - \nabla u_{0} \right) + \int_{\Omega} u_{\varepsilon} \left(u_{\varepsilon} - u_{0} \right)
		&= 
		\int_{\Omega} a_{0} (t,x) \nabla u_{0} \left(  \nabla u_{\varepsilon} - \nabla u_{0} \right) 
		+ \int_{\Omega} u_{0} \left(u_{\varepsilon} - u_{0} \right),\\
		\int_{\Omega} a_{\varepsilon} (t,x) \nabla u_{\varepsilon} \left(  \nabla u_{\varepsilon} - \nabla u_{0} \right) 
		+ \int_{\Omega} \left(u_{\varepsilon}  - u_0 \right)^{2} 
		&= 
		\int_{\Omega} a_{0} (t,x) \nabla u_{0} \left(  \nabla u_{\varepsilon} - \nabla u_{0} \right). 
		\end{align*}
		
		We now subtract $  \int_{\Omega} a_{\varepsilon} (t,x) \nabla u_{0} \left(  \nabla u_{\varepsilon} - \nabla u_{0} \right) $ on both sides, which results
		
		\begin{align*}
		\int_{\Omega} a_{\varepsilon} (t,x)  \left(  \nabla u_{\varepsilon} - \nabla u_{0} \right)^{2} + \int_{\Omega} \left(u_{\varepsilon} - u_{0} \right)^{2} 
		& = 
		\int_{\Omega} [a_{0} (t,x)-a_{\varepsilon}(t,x)] \nabla u_{0} \left(  \nabla u_{\varepsilon} - \nabla u_{0} \right) \\
		& \leq \| a_{\varepsilon} - a_0\|_{\infty}   
		\| \nabla u_{0} \|_{L^{2}} \| \nabla u_{\varepsilon} - \nabla u_{0}\|_{L^{2}} 
		\end{align*}
		
		If $\frac{1}{C}=\min\{ m, 1 \} $, we obtain from the above inequality and using \eqref{uniform_bound},
		
		\begin{align*}
		\|u_{\varepsilon} - u_{0}\|^{2}_{H^{1}}  
		&= 
		\left[ \int_{\Omega}  \left(  \nabla u_{\varepsilon} - \nabla u_{0} \right)^{2} + \int_{\Omega} \left(u_{\varepsilon} - u_{0} \right)^{2} \right]\\
		&\leq C\| a_{\varepsilon}  - a_0  \|_{\infty}  \| \nabla u_{0} \|_{L^{2}}\| \nabla u_{\varepsilon} - \nabla u_{0}\|_{L^{2}}\\
		&\leq \frac12  C^{2} \| a_{\varepsilon} - a_0  \|^{2}_{\infty}  \| \nabla u_{0} \|^{2}_{L^{2}}+ \frac12 \| \nabla u_{\varepsilon} - \nabla u_{0}\|_{L^{2}}^{2}\\ 
		&\leq \frac12  C^{2} \| a_{\varepsilon}- a_0  \|^{2}_{\infty}  \| g \|^{2}_{L^{2}}+ \frac12 \|  u_{\varepsilon} - u_{0}\|_{H^{1}}^{2}.
		\end{align*}
		
		Therefore,
		$
		\|u_{\varepsilon} - u_{0}\|^{2}_{H^{1}}  
		\leq   C \| a_{\varepsilon} - a_0  \|^{2}_{\infty}  \| g \|^{2}_{L^{2}}.
		$
	\end{proof}
\end{lemma}

As an immediate consequence of the previous result, we have the following corollary.

\begin{corollary}\label{op_inv_convergence}
	The operators $A_{\varepsilon}(t)^{-1}: L^{2}(\Omega) \to H^{1}(\Omega)$ are uniformly bounded for $t\in \mathbb{R}$ and $\varepsilon \in [0,1]$ and they converge to $A_0(t)^{-1}$ in the uniform topology. More precisely, for all $\varepsilon\in [0,1]$ and $t\in \mathbb{R}$,
	\begin{align}
	\|A_{\varepsilon}(t)^{-1}\|_{\mathcal{L}  (L^{2}, H^{1} ) } &\leq C, \label{lim_A_e} \\
	\|A_{\varepsilon}(t)^{-1} - A_0(t)^{-1}\|_{\mathcal{L}  (L^{2}, H^{1} ) }
	& \leq C \|a_{\varepsilon}- a_{0}  \|_{ \infty}, \label{conv_A_e}
	\end{align}
	where $C$ does not depend on $\varepsilon$ or $t$.
\end{corollary}

	Inequality \eqref{conv_A_e} is the statement required in \ref{P_5}, with rate of convergence $\eta (\varepsilon)$ given by
	$
	\eta (\varepsilon) := C \| a_{\varepsilon } - a_{0} \|_{\infty}.
	$
	Since \ref{P_1} to \ref{P_5} are satisfied, we conclude that each family of linear operators $\{ A_{\varepsilon}(t), \ t\in \mathbb{R} \}$ generates a linear process $\{U_{\varepsilon} (t,\tau) :L^{2} (\Omega) \to L^{2}(\Omega), \ t\geq\tau, \ \tau \in \mathbb{R} \}$ and from Theorem \ref{T_Conv_Proc}, we obtain that there exist $C,K>0$ such that, for any $\varepsilon\in [0,1]$, $t>\tau$ and $\tau\in \mathbb{R}$, 
	\begin{align*}
	\| U_{\varepsilon}(t,\tau) - U_0(t,\tau) \|_{\mathcal{L}(L^{2})} &\leq C(t-\tau)^{-\theta} e^{K(t-\tau)} \ell (\theta, \varepsilon) ,\\
	\| U_{\varepsilon}(t,\tau) - U_0(t,\tau) \|_{\mathcal{L}(L^{2},H^{1})} &\leq C(t-\tau)^{-\frac12 -\frac{\theta}{2}} e^{K(t-\tau)} \ell(\theta, \varepsilon),
	\end{align*}
	where $\ell (\theta, \varepsilon)= [\| a_{\varepsilon_1}-a_0 \|_{\infty} + \| \nabla_x a_{\varepsilon} - \nabla_x a_0 \|_{\infty}]^{\theta}$ and $\theta \in (0,1)$ is arbitrary, with $C$ depending on the choice of $\theta$.

\subsection{Local well-posedness, global well-posedness and convergence of the solutions}

The results on local and global well-posedness that we present in the sequel can be found in \cite[Section 6]{BCNS_2}. Conditions required for $a_{\varepsilon}(\cdot, \cdot)$ and $f_{\varepsilon}(\cdot, \cdot)$ ensure that problem \eqref{eq_1_abs} admits local solution $u_{\varepsilon}:[\tau,\tau+ T(\varepsilon,\tau, u^{\tau})) \rightarrow H^{1}(\Omega)$ given by
$$	u_{\varepsilon}(t,\tau,u^{\tau}) = U_{\varepsilon}(t,\tau)u^{\tau} + \int_{\tau}^{t} U_{\varepsilon}(t,s) F_{\varepsilon}(s,u_{\varepsilon}(s)) ds,$$					
such that $u_{\varepsilon}(t) \in D=X^{1}$, for all $t\in (\tau, \tau+ T(\varepsilon,\tau, u^{\tau}))$, and $U_{\varepsilon}(\cdot, \cdot): L^{2}(\Omega) \to L^{2}(\Omega)$ is the linear process associated to $\{ A_{\varepsilon}(t), \ t\in \mathbb{R} \}$.

To obtain global well-posedness, we assume that $f_{\varepsilon}$ satisfies a dissipativeness condition: 
\vspace{0.2cm}

\begin{enumerate}[label={\textbf{(D)}},ref=(D)]
	\item \label{D}
	\begin{center} $
		\limsup_{|s| \rightarrow \infty} \left[ sup_{\varepsilon\in [0,1]} \dfrac{f_{\varepsilon}(t,s)}{s} \right] < 1,$
	\end{center}
\end{enumerate}

\vspace{0.2cm}\noindent for all $t\in \mathbb{R}$. The value $1$ comes from the fact that first eigenvalue of $A_{\varepsilon}(t)$ is $\lambda_{\varepsilon,1} (t)=1$. In next lemma, we restate this dissipativeness condition in a manner suitable to applications. Its proof follows directly from the definition of \textit{Limsup}.

\begin{lemma}\label{Lemma_diss_reescrita}
	Suppose that condition \ref{D} holds, then there exists  $\omega_1>0$ such that, for each $\omega \in (0,\omega_1)$, there is a constant $N>0$ such that 
	\begin{equation}\label{est_f_1}
	f_{\varepsilon}(t,s)s \leq (1-\omega) s^{2} +N, \quad \mbox{ for all } s \in \mathbb{R}, \ t\in \mathbb{R}, \ \varepsilon\in [0,1].
	\end{equation}
	
	\noindent Moreover, $N$, $\omega$ and $\omega_1$ are independent of $\varepsilon$.
\end{lemma}

The dissipativeness assumption allows us to obtain global well-posedness, as well as  existence of an absorbing bounded set in $H^{1}(\Omega)$, uniform in $\varepsilon\in [0,1]$.

\begin{theorem}\cite[Theorem 6.13]{BCNS_2}\label{P_A_sets}
	Assume that \ref{A_1} to \ref{A_4} and \ref{D} hold. Let $N, \omega$ be the constants in \eqref{est_f_1} obtained  from the dissipativeness condition \ref{D}. There exists a constant $E>0$ independent of $\varepsilon\in [0,1]$ and of $\tau\in \R$, such that, for any bounded set $B \subset H^{1}(\Omega)$ we can find $T =T(B)>0$, %also independent of $\varepsilon\in [0,1]$ and $\tau \in R$,
	for which
	\[
	\|u_{\varepsilon}(t, \tau, u^{\tau}) \|_{H^{1}} \leq E, \quad \mbox{ for any }u^{\tau} \in B, \ \varepsilon \in [0,1],
	\]
	as long as $t-\tau\geq T$. In particular, the solution of \eqref{eq_1_abs} is globally defined and associated to it there is a nonlinear process $S_{\varepsilon}(t,\tau)$ in $H^{1}(\Omega)$ given by 
	\[
	S_{\varepsilon}(t,\tau) u^{\tau} = u_{\varepsilon}(t,\tau,u^{\tau}) = U_{\varepsilon} (t,\tau) u^{\tau} + \int_{\tau}^{t} U_{\varepsilon}(t,s) F_{\varepsilon}(s,u_{\varepsilon}(s))ds, \mbox{ for all }t\geq \tau.
	\]
\end{theorem}

	Once we proved that the dynamics of all the problems enter a common bounded set $B_{H^{1}} [0,E]$ (the closed ball in $H^{1}(\Omega)$ centered in $0$ and with radius $E$), we can proceed with a cut-off for the nonlinearities $F_{\varepsilon}$, as mentioned in Remark \ref{R_cut-off}. If that is the case, at least close to $B_{H^{1}}[0,E]$, Condition \ref{(NL_1)} holds. As far as Condition \ref{(NL_2)}, we have 
	\begin{align*}
	\| F_{\varepsilon} (t,u ) - F_0 (t,u) \|^{2}_{L^{2}}
	&
	= \int_{\Omega} |f_{\varepsilon} (t,u(x)) - f_{0}(t,u(x)) |^{2} dx \leq \| f_{\varepsilon} - f_0\|_{\infty}^{2} |\Omega|.
	\end{align*}
	
	Therefore, 
	\[
	\sup_{t\in \mathbb{R}} \sup_{u\in H^{1}} \| F_{\varepsilon} (t,u ) - F_0 (t,u) \|_{L^{2}} \leq \|f_{\varepsilon}-f_0 \|_{\infty} |\Omega|^{\frac12} := \gamma(\varepsilon),
	\]
	and Condition \ref{(NL_2)} also holds. We then conclude from Theorem \ref{T_Conv_solution} that the solutions converge as $\varepsilon$ goes to zero, with a rate:
	\begin{equation*}
	\| u_{\varepsilon}(t,\tau, u^{\tau}) - u_0 (t,\tau,u^{\tau}) \|_Y \leq C (t-\tau)^{-\frac{1}{2}-\frac{\theta}{2}} e^{K(t-\tau)} \left[1+ \|u^{\tau}\|_Y\right]\rho(\theta, \varepsilon), 
	\end{equation*}
	
	\vspace{0.1cm} \noindent where $\rho(\theta, \varepsilon) =\max\{  [\| a_{\varepsilon_1}-a_0 \|_{\infty} + \| \nabla_x a_{\varepsilon} - \nabla_x a_0 \|_{\infty}]^{\theta}, \|f_{\varepsilon}-f_0 \|_{\infty} |\Omega|^{\frac12} \}$, $\theta\in (0,1)$ is arbitrary and $C,K>0$ are constants independent of $\varepsilon \in [0,1], \ \tau \in \R, \ t >\tau$ or $u^{\tau}\in Y$, but dependent on the choice of $\theta.$

\subsection{Existence of pullback attractor and its upper-semicontinuity}

The existence of pullback attractor was also obtained in \cite[Theorem 6.14]{BCNS_2} for a single problem, rather than a family of problems. However, in the proof of this theorem, the authors were able to find a compact absorbing set in $H^{1}(\Omega)$ depending only on the constants $\omega$, $N$ obtained in Lemma \ref{Lemma_diss_reescrita} and the growth $\rho$ for the nonlinearity. Since they are all uniform for $\varepsilon \in [0,1]$, we can state the next theorem as a consequence of the result in \cite[Theorem 6.14]{BCNS_2}.

\begin{theorem}
	\label{T_existence_attractor}
	Assume that \ref{A_1} to \ref{A_4} and \ref{D} hold. Let $N, \omega$ be the constants in \eqref{est_f_1} obtained  from the dissipativeness condition \ref{D}. Then the nonlinear process $S_{\varepsilon}(t,\tau) = u_{\varepsilon}(t,\tau,\cdot)$ in $H^{1}(\Omega)$  has a pullback attractor $  \{\mathcal{A}_{\varepsilon} (t), \  t\in \mathbb{R} \}$ in $H^{1}(\Omega)$. Moreover, there exists a compact set $K \subset H^{1}(\Omega)$ such that 
	\begin{equation}\label{Unif_comp}
	\left[\scalebox{1}{$\bigcup_{\varepsilon \in [0,1]}  \bigcup_{t\in \mathbb{R}}$} \mathcal{A}_{\varepsilon} (t) \right] \subset K.
	\end{equation} 
\end{theorem}	
	
	\
	
	 From Corollary \ref{C_NL_process_convergence}, we obtain continuity of the family $\{S_{\varepsilon}(\cdot,\cdot)\}_{\varepsilon \in [0,1]}$ and from \eqref{Unif_comp}, we conclude that 
	\[
	\overline{\left[\scalebox{1}{$\bigcup_{\varepsilon \in [0,1]}  \bigcup_{t\in \mathbb{R}}$} \mathcal{A}_{\varepsilon} (t) \right]}
	\]
	is relatively compact. Those are the conditions in \cite[Theorem 3.6]{Carvalho} necessary to ensure upper-semicontinuity of the family $\{\mathcal{A}_{\varepsilon} (t), t\in \mathbb{R}\}_{\varepsilon\in [0,1]}$ at $\varepsilon=0$.

	\begin{corollary}
		Under conditions of Theorem \ref{T_existence_attractor}, the family of pullback attractos $\{\mathcal{A}_{\varepsilon} (t), t\in \mathbb{R}\}_{\varepsilon\in [0,1]}$ is upper-semicontinuous at $\varepsilon=0$.
	\end{corollary}

\vspace{0.1cm}
	
%%%%%%%%%%%%%%%%%%%%%%%%%%%%%%%%%%%%%%%%
%SECTION 5 - PERTURBATION OF FRACTIONAL POWER
%%%%%%%%%%%%%%%%%%%%%%%%%%%%%%%%%%%%%%%%%%%%%%%%%%%%%

\section{Application to a nonautonomous strongly damped wave equations and its fractional approximations}\label{S_application_II}

As a second application, we consider the nonautonomous strongly damped wave equation subjected to Dirichlet boundary conditions
\begin{equation}\label{SDWE}
\begin{split}
& u_{tt}+(-a(t) \Delta_D) u + 2 (-a(t)\Delta_D)^{\frac12} u_t = f(t,u), \hspace{1cm} x\in \Omega, \ t>\tau, \\
& u(t,x) = 0, \hspace{6.4cm}   x\in \partial \Omega, \ t\geq \tau, \\
& u(\tau,x) = u^{\tau}(x), \ u_t (\tau, x)=v^{\tau}(x), \hspace{2.98cm} x\in \overline{\Omega}, \ \tau \in \mathbb{R},
\end{split}
\end{equation}
where $\Omega \subset \mathbb{R}^{n}$, $n\geq 3$, is a bounded smooth domain, $\Delta_D$ is the Laplacian operator with Dirichlet boundary condition and $f: \R \times \R \to \R$ a nonlinearity. We shall assume the following additional condition:

\vspace{0.2cm}

\begin{enumerate}[label={(B)},ref=(B)]

	\item \label{B_1} The function $a :\R \to  \R^{+}$ is positive and has its image in a bounded interval of the form $[a_0,a_1] \subset (0,\infty)$. We also assume that it is H\"{o}lder continuous with an exponent $\delta\in (0,1]$, that is, there exists a constant $C>0$ such that
	\begin{equation*}\label{HC_a}
	|a(t) - a(s)| \leq C|t-s|^{\delta}, \quad \mbox{for all } t,s \in \mathbb{R}.
	\end{equation*}

\end{enumerate}

\vspace{0.2cm}

	Let $E=L^{2}(\Omega)$ and denote by $A(t):D(A(t)) \subset E \to E$ the operator
	\begin{equation}\label{A(t)}
	A(t)u = -a(t)\Delta_{D}u, \quad \mbox{for} \quad u\in D(A(t)) = D (\Delta_D)  = H^{2}(\Omega)\cap H^{1}_0(\Omega),
	\end{equation}
	where $D(\Delta_{D})$ stands for the domain of the Laplacian with Dirichlet boundary conditions. As expected, the multiplication by a real-function $a(t)$ does not change the domain of the Laplacian.
	
	Therefore, this linear operator $A(t)$ has a time-independent domain and from the well-known properties of the Laplacian operator \cite{pazy} and the fact that $a(t) \geq a_0,$ for all $t\in \mathbb{R}$, we deduce that $A(t)$ is a positive operator, self-adjoint, sectorial and $-A(t)$ generates a compact analytic $C_0-$ semigroup in $E$.

	Consequently, fractional powers of $A(t)$ in the sense of Amman \cite{amann} are well-defined. We shall denote by $A(t)^{\alpha}$ the power of the linear operator $A(t)$. One can easily deduce from the expressions for fractional power of linear operators that 
	\[
	A(t)^{\alpha} = (-a(t)\Delta_D)^{\alpha} = [a(t)]^{\alpha}  (-\Delta_D)^{\alpha}, \quad \mbox{ for all } t\in \R \mbox{ and } \alpha \in (0,1],
	\]
	and the domain of $A(t)^{\alpha}$ is the same as the domain of $(-\Delta_D)^{\alpha}$, that is,
	\[
	D(A(t)^{\alpha}) = D((-\Delta)^{\alpha}), \quad \mbox{for all } t\in \R.
	\]
	
	We then define a scale of Banach spaces given by the fractional powers  $(-\Delta_D)^{\alpha}$, $\alpha \in (0,1]$,
	\[
	E^{\alpha} = (-\Delta_D)^{\alpha} \quad \mbox{eqquiped with the norm } \|\cdot \|_{E^{\alpha}} = \| (-\Delta_D)^{\alpha} \cdot \|_{L^{2}}.
	\]
\noindent In particular, $E^{0} = E = L^{2}(\Omega)$, $E^{\frac12} = H_0^{1}(\Omega)$ and $E^{1} = H^{2}(\Omega) \cap H_0^{1}(\Omega)$ (see \cite{CholewaDlotko}Theorem 3.6). From the boundedness of $a(\cdot)$, there exists $0<m<M$ such that $[a(t)]^{\alpha} \in [m,M]$ for all $t\in \mathbb{R}$ and $\alpha \in (0,1]$. Therefore, 
	\begin{equation*}\label{eq_norms}
	m \|u \|_{E^{\alpha}} \leq \|A(t)^{\alpha} u \|_{L^{2}} \leq M \|u \|_{E^{\alpha}}
	\end{equation*}
	and the graph norm associated to the linear operator $A(t)^{\alpha}$ is equivalent to the norm in $E^{\alpha}$. With the above set up and  taking $u_t =v$, Problem \eqref{SDWE} can be written in the following abstract form

	\begin{equation}\label{SDWE_abstract}
	\frac{d}{dt} \begin{bmatrix}
	u \\
	v
	\end{bmatrix}
	+
	\Lambda(t)
	\begin{bmatrix}
	u \\
	v
	\end{bmatrix}
	=
	F\left(t,
	\scalebox{0.9}{$
	\begin{bmatrix}
	u \\
	v
	\end{bmatrix}$}
	\right), 
	\ t>\tau ; \qquad
	%
	%
	%INITIAL CONDITION
	%
	%
	\begin{bmatrix}
	u(\tau) \\
	v(\tau)
	\end{bmatrix} 
	=
	\begin{bmatrix}
	u^{\tau}\\
	v^{\tau}
	\end{bmatrix} \in E^{\frac12} \times E ,
	\end{equation} 
	where $\Lambda(t) : D(\Lambda(t)) \subset E^{\frac12} \times E  \to E^{\frac12} \times E$ is the linear operator defined in $D(\Lambda(t)) = D = E^{1} \times E^{\frac12}$ and given by
	\begin{equation}\label{Lambda(t)}
	\Lambda(t)
	\begin{bmatrix}
	u\\
	v
	\end{bmatrix}=
	\begin{bmatrix}
	0 & -I\\
	A(t) & 2 A(t)^{\frac12}
	\end{bmatrix}
	\begin{bmatrix}
	u\\v
	\end{bmatrix}
	=
	\begin{bmatrix}
	-v\\
	A(t)u + 2A(t)^{\frac12}v
	\end{bmatrix},
	\end{equation}
	and $F$ is the nonlinearity given by
	\begin{equation*}
	F\left(t,
\scalebox{0.9}{$
	\begin{bmatrix}
	u \\
	v
	\end{bmatrix}$}
	\right)
	=
	\begin{bmatrix}
	0 \\
	f(t,u)
	\end{bmatrix}.
	\end{equation*}

	We have the following result proved in \cite[Lemma 8.1]{BezerraCarvalhoNascimento2020} concerning spectral properties of $\Lambda(t)$ and the calculus of its fractional powers.
	
	\begin{proposition}\label{L_prop_Lambda}
	If $A(t)$ and $\Lambda(t)$ are as in \eqref{A(t)} and in \eqref{Lambda(t)}, respectively, then the following properties hold:
	
	\begin{enumerate}
		\item $\Lambda(t)$ is a positive operator with fractional powers denoted by $\Lambda(t)^{\alpha}$, $\alpha\in (0,1]$. 
		
		\item There exists $C>0$ and $\varphi \in (\frac{\pi}{2}, \pi)$ (independent of $\alpha$) such that, for any $\alpha \in (0,1]$ and $t\in \mathbb{R}$, 
		\[
		\Sigma_{ \varphi} \cup\{0\} \subset \rho(-\Lambda(t)^{\alpha}),
		\]
		and the following estimates hold
		\[
		\| (\lambda +\Lambda(t)^{\alpha})^{-1}\|_{\mathcal{L}(E^{\frac12}\times E)} \leq \frac{C}{1+|\lambda|},  \quad \mbox{ for all }\lambda \in \Sigma_{ \varphi} \cup\{0\}.
		\]
		
		\noindent Therefore, each $\Lambda(t)^{\alpha}$ is sectorial in $E^{\frac12} \times E$ and $-\Lambda(t)^{\alpha}$ generates an analytic semigroup in $\mathcal{L}( E^{\frac12} \times E)$.

		\item Given any $\alpha \in (0,1]$, we have the following explicitly expression for the fractional powers of $\Lambda(t)$:
		
		\begin{equation}\label{Frac_expressions}
		\Lambda(t)^{\alpha} =
		\begin{bmatrix}
		(1-\alpha) A(t)^{\frac{\alpha}{2}} & -\alpha A(t)^{\frac{-1+\alpha}{2}}\\
		\alpha A(t)^{\frac{1+\alpha}{2}} & (1+\alpha) A(t)^{\frac{\alpha}{2}}
		\end{bmatrix}
		\quad
		\mbox{and}
		\quad 
		\Lambda(t)^{-\alpha} =
		\begin{bmatrix}
		(1+\alpha) A(t)^{-\frac{\alpha}{2}} & \alpha A(t)^{\frac{-1-\alpha}{2}}\\
		-\alpha A(t)^{\frac{1-\alpha}{2}} & (1-\alpha) A(t)^{-\frac{\alpha}{2}}
		\end{bmatrix}
		\end{equation}
	\end{enumerate} 
	
	\end{proposition}
	
	\

	We shall consider fractional versions of Problem \eqref{SDWE_abstract}, given by 
	\begin{equation}\label{Fractional_SDWE_abstract}
	\frac{d}{dt} \begin{bmatrix}
	u_{\alpha} \\
	v_{\alpha}
	\end{bmatrix}
	+
	\Lambda(t)^{\alpha}
	\begin{bmatrix}
	u_{\alpha} \\
	v_{\alpha}
	\end{bmatrix}
	=
	F_{\alpha}\left(t,
	\scalebox{0.9}{$
		\begin{bmatrix}
		u_{\alpha} \\
		v_{\alpha}
		\end{bmatrix}$}
	\right), 
	\ t>\tau ; \qquad
	%
	%
	%INITIAL CONDITION
	%
	%
	\begin{bmatrix}
	u_{\alpha}(\tau) \\
	v_{\alpha}(\tau)
	\end{bmatrix} 
	=
	\begin{bmatrix}
	u^{\tau}\\
	v^{\tau}
	\end{bmatrix} \in E^{\frac12} \times E ,
	\end{equation} 
	where 
	$\alpha \in (0,1]$, $\Lambda(t)^{\alpha}$ is the fractional power of $\Lambda(t)$ and 
	\[
	F_{\alpha}\left(t,
	\scalebox{0.9}{$
		\begin{bmatrix}
		u_{\alpha} \\
		v_{\alpha}
		\end{bmatrix}$}
	\right) =
	\begin{bmatrix}
	0 \\
	f_{\alpha}(t,u_{\alpha})
	\end{bmatrix}.
	\]
	
	By analyzing Expression \eqref{Frac_expressions} for the linear operator $\Lambda(t)^{\alpha}$, we see that  $\alpha= 1$ recovers the original expression for $\Lambda(t)$, so we might expect that as we make $\alpha \to 1^{-}$, the fractional problem \eqref{Fractional_SDWE_abstract} approaches \eqref{SDWE_abstract} in a certain sense. We shall verify that this is the case, that is, we prove in the sequel that conditions \ref{P_1} to \ref{P_5} hold for this family (in $\alpha$) of problems. Therefore, if $U_{\alpha}(t,\tau)$ is the linear process associated to $\{\Lambda(t)^{\alpha}, t \in \R \}$, then we  obtain its convergence to the linear process $U(t,\tau)$ associated to $\{\Lambda (t), t\in \R \}$.

	A slightly different singularly nonautonomous wave equation  and its fractional perturbations  were also considered in  \cite{BezerraNascimento_2019}. However, the lack of results on the linear process convergence prevented the authors  to proceed with the analysis beyond spectrum convergence of the fractional operators. They were not able to obtain convergence of the associated linear processes or of the solutions, as we shall do in the sequel.
	 
	Using the notation developed in Section \ref{S:Setting_and_main_results}, we will consider $Y=X= E^{\frac12} \times E$. Conditions \ref{P_1} and \ref{P_2} follows directly from Proposition \ref{L_prop_Lambda}. For condition \ref{P_3}, we shall need the following technical lemma.

	\begin{lemma}\label{L_HC_a(.)}
	Let $a(\cdot)$ be the function satisfying \ref{B_1}. For each $\omega > 0$, the functions $[a(\cdot)]^{\omega}$ and $[a(\cdot)]^{-\omega}$ are also H\"{o}lder continuous with H\"{o}lder exponent $\delta$, that is, for all $t,s\in \R$,
	\[
	|[a(t)]^{\omega}-[a(s)]^{\omega}| \leq C_1 |t-s|^{\delta} \quad \mbox{ and } \quad |[a(t)]^{-\omega}-[a(s)]^{-\omega}| \leq C_2 |t-s|^{\delta}.
	\]
	\begin{proof}
		Let $\phi: \R^{+} \to \R^{+}$ be given by $\phi(s)=s^{\omega}$. From the mean value theorem and the fact that $a(t) \in [a_0,a_1]$, for all $t\in \R$, we obtain, for some $\theta$ between $a(s)$ and $a(t)$, in particular $\theta \in  [a_0,a_1]$, that
		\[
		\begin{split}
		|[a(t)]^{\omega} - [a(s)]^{\omega}|=	|\phi(a(t)) - \phi(a(s))| \leq |\phi'(\theta)| |a(t)-a(s)| \leq \theta^{\omega-1}|a(t)-a(s)| 
		\leq C_1 |t-s|^{\delta}
		\end{split},
		\]
		for all $t,s \in \R$. Moreover, 
		\[
		|[a(t)]^{-\omega} - [a(s)]^{-\omega} |= \left| \frac{  [a(s)]^{\omega} - [a(t)]^{\omega} }{[a(t)]^{\omega}[a(s)]^{\omega}   } \right|
		\leq C|[a(s)]^{\omega} - [a(t)]^{\omega}| \leq C_2|t-s|^{\delta}, \quad \mbox{ for all }t,s\in \R.
		\]
	
	\end{proof}
	\end{lemma}

	Condition \ref{P_3} can now be verified.

	\begin{lemma}
	Assume that \ref{B_1} holds and let $\Lambda(t)$ be the linear operator in \eqref{Lambda(t)}. If $\Lambda(t)^{\alpha}$ denotes the fractional powers of $\Lambda(t)$, then $\mathbb{R} \ni t \to  \Lambda(t)^{\alpha} \Lambda(\tau)^{-\alpha} \in \mathcal{L}(E^{\frac12}\times E)$ is $\delta-$H\"{o}lder continuous, {uniformly in $\alpha\in (0,1]$} and $\tau \in \R$. In other words, there exists $C>0$ such that 
		\[
		\|   [\Lambda(t)^{\alpha} - \Lambda(s)^{\alpha}] \Lambda(\tau)^{-\alpha} \|_{ \mathcal{L} (E^{\frac12}\times E)   }
		\leq 
		C|t-s|^{\delta}, \quad \mbox{ for all } t,  s,  \tau \in \mathbb{R} \mbox{ and  } \alpha \in (0,1].
		\]

	\begin{proof}
	Applying expression \eqref{Frac_expressions} for the fractional powers, we deduce that 
	\[
	 [\Lambda(t)^{\alpha} - \Lambda(s)^{\alpha}] \Lambda(\tau)^{-\alpha} = 
	 \begin{bmatrix}
	 \Theta_{11} & \Theta_{11}\\
	 \Theta_{21} & \Theta_{22}
	 \end{bmatrix},
	\]
	where 
	\begin{align*}
	\Theta_{11} &= \alpha^{2} \left( [a(t)]^{\frac{\alpha-1}{2}} - [a(s)]^{\frac{\alpha-1}{2}} \right) [a(\tau)]^{\frac{1-\alpha}{2}} 
	+
	 (1-\alpha^{2})
	 \left( [a(t)]^{\frac{\alpha}{2}} - [a(s)]^{\frac{\alpha}{2}} \right) [a(\tau)]^{\frac{-\alpha}{2}}, \\
	 \Theta_{12} & = \alpha (1-\alpha) (-\Delta_D)^{-\frac12}
	 \left\{ 
	 \left( [a(t)]^{\frac{\alpha}{2}} - [a(s)]^{\frac{\alpha}{2}} \right) [a(\tau)]^{\frac{-1-\alpha}{2}} 
	 +
	 \left( [a(s)]^{\frac{\alpha-1}{2}} - [a(t)]^{\frac{\alpha-1}{2}} \right) [a(\tau)]^{\frac{-\alpha}{2}} 
	 \right\},\\
	  \Theta_{21} & = - \alpha (1+\alpha) (-\Delta_D)^{\frac12}
	 \left\{ 
	 \left( [a(t)]^{\frac{\alpha}{2}} - [a(s)]^{\frac{\alpha}{2}} \right) [a(\tau)]^{\frac{1-\alpha}{2}} 
	 +
	 \left( [a(s)]^{\frac{\alpha+1}{2}} - [a(t)]^{\frac{\alpha+1}{2}} \right) [a(\tau)]^{\frac{-\alpha}{2}} 
	 \right\},\\
	 \Theta_{22} &= \alpha^{2} \left( [a(t)]^{\frac{\alpha+1}{2}} - [a(s)]^{\frac{\alpha+1}{2}} \right) [a(\tau)]^{\frac{-\alpha-1}{2}} 
	 +
	 (\alpha+1)(\alpha-1)
	 \left( [a(t)]^{\frac{\alpha}{2}} - [a(s)]^{\frac{\alpha}{2}} \right) [a(\tau)]^{\frac{-\alpha}{2}}.
	\end{align*}
	
	We must obtain H\"{o}lder continuity of each entry in its appropriate space, that is, $\Theta_{11}$ in $\mathcal{L}(E^{\frac12})$, $\Theta_{12}$ in $\mathcal{L}(E, E^{\frac12})$, $\Theta_{21}$ in $\mathcal{L}(E^{\frac12}, E)$ and $\Theta_{22}$ in $\mathcal{L}(E)$. All of them are similar and follows from Lemma \ref{L_HC_a(.)}. We illustrate how to proceed with $\Theta_{21}$.
	\begin{align*}
	 \| \Theta_{21} \|_{\mathcal{L} (E^{\frac12},E)  }
	 &\leq 
	 \| \Theta_{21} (-\Delta_{D})^{-\frac12} \|_{\mathcal{L}(E)} 
	 \leq 
	 \alpha(1+\alpha) 
	\left(
	| [a(t)]^{\frac{\alpha}{2}} - [a(s)]^{\frac{\alpha}{2}} | a_1^{\frac{1-\alpha}{2}} 
	 +
	 |[a(s)]^{\frac{\alpha+1}{2}} - [a(t)]^{\frac{\alpha+1}{2}} | a_0^{\frac{-\alpha}{2}} 
	\right)\\
	& \leq C_{21} |t-s|^{\delta}.
	\end{align*}
	
	The other entries follow analogously.
	\end{proof}
	\end{lemma}

	In a similar way, from the Expressions \eqref{Frac_expressions} for the fractional powers, we can deduce property \ref{P_4} for the family of linear operators $\{\Lambda(t)^{\alpha}, \ t\in \mathbb{R} \}_{\alpha \in (0,1]}$.

	\begin{lemma}
	Assume that \ref{B_1} holds and let $\Lambda(t)$ be the linear operator in \eqref{Lambda(t)} and $\Lambda(t)^{\alpha}$ its fractional powers. There exists a continuous function $\xi: (0,1]\to \R^{+}$ with $\xi(1)=0$ such that 
	\begin{equation}\label{P_4_ex_2}
	\sup_{t\in \R} \| \Lambda(t)^{\alpha}\Lambda(\tau)^{-\alpha} - \Lambda(t) \Lambda(\tau)^{-1} \|_{\mathcal{L} (E^{\frac12}\times E) } \leq \xi (\alpha).
	\end{equation}
	\begin{proof}
	Applying expressions \eqref{Frac_expressions} for the fractional powers, we deduce 
	\[
	\Lambda(t)^{\alpha}\Lambda(\tau)^{-\alpha} - \Lambda(t) \Lambda(\tau)^{-1} 
	=
	\begin{bmatrix}
	\Gamma_{11}&\Gamma{12}\\
	\Gamma{21}& \Gamma{22}
	\end{bmatrix},
	\]
	where
	\begin{align*}
	&\Gamma_{11} = (1-\alpha^{2}) [a(t)]^{\frac{\alpha}{2}}   [a(\tau)]^{\frac{\alpha-1}{2}}  
	+ 
	(\alpha^{2}-1)  [a(t)]^{\frac{\alpha-1}{2}} [a(\tau)]^{\frac{\alpha}{2}}
	+
	\left(  [a( \tau )]^{\frac{ \alpha}{2 }}  [a( t )]^{\frac{\alpha-1 }{2 }}  -  [a( \tau )]^{\frac{ 1}{ 2}} \right),
	\\
	&\Gamma_{12}= (1-\alpha)\alpha
	 [a(t)]^{\frac{\alpha-1}{2}} [a( \tau )]^{\frac{ -\alpha-1}{2}} 
	 \left(  [a(t )]^{\frac{ 1}{2}}  -   [a( \tau )]^{\frac{ 1}{2 }}  \right) (-\Delta_D)^{-\frac12},
	\\
	&\Gamma_{21}=\alpha (1+\alpha) (-\Delta_{D})^{\frac12}
	\left(
	 a( t ) \left\{   [a( t )]^{\frac{\alpha-1 }{ 2}}  - \frac{2}{\alpha(1+\alpha)} \right\}
	 +
	 \left\{
	 \frac{2}{\alpha(1+\alpha)} -  [a( \tau )]^{\frac{1-\alpha }{2 }}
	 \right\}
	\right),
	\\
	&\Gamma_{22}=  [a( t )]^{\frac{ \alpha}{ 2}}  [a(  \tau )]^{\frac{ -1-\alpha}{2 }} 
	\left(
	(1+\alpha^{2}) [a( t)]^{\frac{1}{2}} - \alpha^{2}  [a(\tau )]^{\frac{1 }{2 }}
	\right)
	-
	 [a( t)]^{\frac{ 1}{2 }} [a( \tau )]^{-1}
	 \left(
	 2 [a( t)]^{\frac{1 }{2 }}- [a( \tau )]^{\frac{ 1}{2 }}
	 \right).
	\end{align*}
	
	Note that each entry goes to zero continuously as $\alpha \to 1^{-}$. The presence of $(-\Delta_D)^{-\frac12}$ in $\Gamma_{12}$ or $(-\Delta_D)^{\frac12}$ in $\Gamma_{21}$ do not represent any problem in the estimates. Actually, since $\Gamma_{12}$ is estimated in $\mathcal{L}(E, E^{\frac12})$ and $\Gamma_{21}$  in $\mathcal{L} (E^{\frac12}, E)$ we have
	\[
	\| \Gamma_{12} \|_{\mathcal{L} (E,E^{\frac12})  } = \| (-\Delta_D)^{\frac12}\Gamma_{12}\|_{\mathcal{L}(E)} 
	\quad \mbox{and }\quad
	\| \Gamma_{21} \|_{\mathcal{L} (E^{\frac12},E)  } = \| \Gamma_{12} (-\Delta_D)^{-\frac12}\|_{\mathcal{L}(E)} 
	\]
	and those powers of the Laplacian disappear when we estimate those terms. Therefore, there exists a continuous function $\xi:(0,1]\to \R^{+}$ with $\xi(1)=0$ such that \eqref{P_4_ex_2} holds.
	\end{proof}
	\end{lemma}

	Even though we are able to prove that \ref{P_4} holds, we cannot obtain an explicit formulation for $\xi$ in \eqref{P_4_ex_2}, since it depends on the expression of $a(\cdot)$. Lastly, Condition \ref{P_5} holds following the same proof of \cite[Theorem 3.1]{BezerraCarvalhoNascimento2020}.
	
	\begin{lemma}\cite[Theorem 3.1]{BezerraCarvalhoNascimento2020}
		Let $\Lambda(t)$ be the linear operator in \eqref{Lambda(t)} and $\Lambda(t)^{\alpha}$, $\alpha\in (0,1]$ its fractional. There exists a constant $C>0$, independent of $\alpha$ and $t\in \R$, such that
		\[
		\| \Lambda(t)^{-\alpha} - \Lambda(t)^{-1} \|_{\mathcal{L}(E^{\frac12} \times E)} \leq C(1-\alpha). 
		\]
	\end{lemma}

	Since conditions \ref{P_1} to \ref{P_5} hold, we have the following result, which is a restatement of Theorem \ref{T_Conv_Proc}.

	\begin{theorem}
		Let $\{U(t,\tau) \in \mathcal{L}(E^{\frac12}  \times E ), \ t\geq \tau\}$ be the linear process associated to $\{\Lambda(t), \ t\in \R\}$ and $\{U_{\alpha}(t,\tau)\in \mathcal{L}(E^{\frac12}  \times E ), \ t\geq \tau\}$ the linear process associated to $\{\Lambda(t)^{\alpha}, \ t\in \R\}$, ${\alpha\in (0,1]}$. For any $\theta \in (0,1)$, there exists constants $K,C>0$, independent of $\alpha \in {(0,1]}$, $\tau \in \R$ and $t>\tau$, such that 
		\[
		\| U_{\alpha}(t,\tau) - U(t,\tau) \|_{\mathcal{L} (E^{\frac12} \times E) } \leq C(t-\tau)^{-\theta} e^{K(t-\tau)} \ell (\theta, \alpha),
		\]
		where $\ell (\theta, \alpha) = \max\{ (1-\alpha)^{\theta}, [\xi(\alpha)]^{\theta}\}$.
	\end{theorem}	

	Under assumptions on boundedness and Lipschitz continuity for the family of nonlinearities $f_{\alpha} : \R \times R \to \R$, $\alpha \in (0,1]$, as well as some convergence assumption of $f_{\alpha}$ to $f$ as $\alpha \to 1^{-}$,  we derive conditions \ref{(NL_1)} and \ref{(NL_2)} for $F_{\alpha}: E^{\frac12} \times E \to E^{\frac12} \times E$, as we did in Section \ref{S:Case_I_application_I}. Then we could obtain convergence of the solutions $u_{\alpha}(t,\tau, [u^{\tau}, v^{\tau}])$ of Problem \eqref{Fractional_SDWE_abstract} to the solution $u_{}(t,\tau, [u^{\tau}, v^{\tau}])$ of Problem \eqref{SDWE_abstract} in $E^{\frac12} \times E$, as $\alpha \to 1^{-}$ as a consequence of Theorem \ref{T_Conv_solution}.

\section{Disclosures and declaration}

The author declares no conflicts of interest. Maykel Boldrin Belluzi has received a research grant from FAPESP, Brazil, process number 2022/01439-5.

Moreover, data sharing is not applicable to this article as no new data were created or analyzed in this study.
	
%%%%%%%%%%%%%%%%%%%%%%%%%%%%%%%%%%%
									%BIBLIOGRAFIA
%%%%%%%%%%%%%%%%%%%%%%%%%%%%%%%%%%%
	
% acm é um comum usado pela matemática.
%
%Ver outros padrões de citação em:  https://www.reed.edu/cis/help/LaTeX/bibtexstyles.html

%\bibliographystyle{acm}\bibliography{reference}

\end{document}